%% file: paper.tex
\documentclass{article}

\usepackage[utf8]{inputenc}
\usepackage[T1]{fontenc}

\usepackage{arxiv}
\renewcommand{\headeright}{}
\renewcommand{\undertitle}{}
\renewcommand{\shorttitle}{Benchmarks in Leipzig}
\usepackage{amsmath, amssymb}
\usepackage{natbib}
\usepackage{hyperref}
\usepackage{booktabs}
\usepackage{multirow}
\usepackage{caption}
\usepackage{xcolor}
\usepackage{microtype}
\usepackage{enumitem}
\setlist[itemize]{leftmargin=1em}
\usepackage[most]{tcolorbox}
\usepackage{cleveref}
\usepackage{doi}
\usepackage{xspace}
\usepackage{array}

\usepackage{tocloft}
\newcommand{\gpt}{\texttt{GPT-5.5}\xspace}
\newcommand{\gemini}{\texttt{Gemini 3.1 Pro}\xspace}
\newcommand{\grok}{\texttt{Grok 4.3}\xspace}
\newcommand{\deepseek}{\texttt{DeepSeek-V4-Pro}\xspace}
\newcommand{\claude}{\texttt{Claude Opus 4.7}\xspace}

\newcommand{\astra}{\texttt{GPT-6 Astra}\xspace}
\newcommand{\astrapro}{\texttt{GPT-6 Astra Pro}\xspace}
\newcommand{\sol}{\texttt{GPT-5.6 Sol}\xspace}
\newcommand{\solpro}{\texttt{GPT-5.6 Sol Pro}\xspace}
\newcommand{\opusfive}{\texttt{Claude Opus 5}\xspace}
\newcommand{\fable}{\texttt{Claude Fable 5.1}\xspace}
\newcommand{\gemflash}{\texttt{Gemini 3.8 Flash}\xspace}
\newcommand{\gptsixsol}{\texttt{GPT-6 Sol}\xspace}
\newcommand{\opusfivefive}{\texttt{Claude Opus 5.5}\xspace}

\newtcolorbox{problembox}[1][]{
  enhanced,
  breakable,
  colframe = gray!30,
  colback  = gray!5,
  coltitle = gray!20!black,
  fonttitle = \bfseries,
  title    = {#1},
  boxrule  = 0.6pt,
  arc      = 1pt,
  left     = 8pt,
  right    = 8pt,
  top      = 4pt,
  bottom   = 4pt,
}

\newtcolorbox{findingsbox}[1][Key findings]{
  enhanced,
  breakable,
  colframe = gray!30,
  colback  = gray!5,
  coltitle = gray!20!black,
  fonttitle = \bfseries,
  title    = {#1},
  boxrule  = 0.6pt,
  arc      = 1pt,
  left     = 8pt,
  right    = 8pt,
  top      = 4pt,
  bottom   = 4pt,
}

\hypersetup{
  colorlinks = true,
  linkcolor  = blue!50!black,
  urlcolor   = blue!50!black,
  citecolor  = blue!50!black,
}

\title{Benchmarks in Leipzig\\[0.4em]
       \large A collection of questions in research-level mathematics}

\author{%
\parbox{\dimexpr\textwidth-6pc\relax}{\normalfont\raggedright
{\large\bfseries Project Supervisor}\\[0.4ex]
Christian Stump$^{1}$\thanks{Corresponding author: \href{mailto:christian.stump@rub.de}{christian.stump@rub.de}. For questions or model evaluation requests: \href{mailto:christian@sciencebench.ai}{christian@sciencebench.ai}.}\\[1.5ex]
{\large\bfseries Event Organizers}\\[0.4ex]
Veronica Calvo Cortes$^{2}$, Christian Stump$^{1}$, Bernd Sturmfels$^{2}$\\[1.5ex]
{\large\bfseries Benchmark Contributors}\\[0.4ex]
Andrei~Balakin$^{3}$, Mikl\'os B\'ona$^{4}$, Marie-Charlotte~Brandenburg$^{1}$, Clara~Briand$^{5}$, Veronica~Calvo~Cortes$^{2}$, Shelby~Cox$^{2}$, Jesus~A.~De~Loera$^{6}$, Danai~Deligeorgaki$^{7}$, Hannah~Friedman$^{8}$, Tim~Gehrunger$^{9}$, Chiara~Giardino$^{1}$, Stephen~Griffeth$^{10}$, Baran~Hashemi$^{2}$, Elena~Hoster$^{1}$, Alexander~Ivanov$^{1}$, Nupur~Jain$^{1}$, Aryaman~Jal$^{1}$, Leonie~Kayser$^{2}$, Joris~Koefler$^{2}$, Kevin~Kühn$^{3}$, Mario~Kummer$^{11}$, Felix~Lotter$^{2}$, René~Marczinzik$^{12}$, Victor~S.~Miller, Alejandro~Morales$^{13}$, Greta~Panova$^{14}$, Gianni~Petrella$^{15}$, Nathan~Pflueger$^{16}$, Lakshmi~Ramesh$^{17}$, Nikolas~Rieke$^{18}$, Carlos~Rodriguez$^{2}$, Andrea~Rosana$^{2}$, Flavio~Salizzoni$^{2}$, Otto~T.P.~Schmidt$^{2}$, Sven~Ulf~Schmitz$^{12}$, Lina~Maria Simbaqueba~Marin$^{19}$, Luca~Sodomaco$^{2}$, Christian~Stump$^{1}$, Bernd~Sturmfels$^{2}$, Alexander~Taveira~Blomenhofer$^{20}$, Simon~Telen$^{2}$, Philipp~Tuchel$^{1}$, Emil~Verkama$^{21}$, Carl~Felix~Waller$^{22}$, Julian~Weigert$^{2,23}$, Annette~Werner$^{24}$, Nathan~Williams$^{25}$, and Claudius~Zibrowius$^{1}$%
}%
}

\date{Version~2, \today}

\begin{document}

\maketitle

\begin{abstract}
Between April 1 and May 15, 2026, a group of 49 mathematicians compiled a dataset of research-level mathematics questions with known answers.
Most of the work was done during the three-day workshop \textbf{Benchmarks in Leipzig} with~35 participants at the Max Planck Institute for Mathematics in the Sciences in Leipzig, Germany.
We present the resulting collection of 100~questions.
We evaluated these questions in three stages: a single attempt by five state-of-the-art LLMs and their predecessors, followed by a~20-runs-per-model evaluation with three of these models, and finally a 3-run attempt with two heavy-thinking models.
After Stage~1,~41 questions remained completely unsolved; after Stage~2, this count dropped to~16; and we concluded Stage~3 with only~2 unsolved questions.
This demonstrates that the mathematical reasoning capabilities of LLMs are becoming impressive.
In September 2026, we added a fourth stage in which the next generation of models attempted all~100 questions once more, after which only~1 question remains unsolved.
\end{abstract}

\setcounter{tocdepth}{2}
\tableofcontents
\medskip

\section*{Acknowledgments}

We thank the Max Planck Institute for Mathematics in the Sciences in Leipzig, Germany, for its hospitality, and for providing financial support for travel and accommodation for external participants.
The institute also funded the API calls to the AI models for the benchmark questions as well as for the chat system.
Contributors did not receive any payment for the questions they submitted.
We also thank \href{https://surgehq.ai/}{Surge AI} for performing external model runs on the benchmark questions.

\section{Introduction}
\label{sec:intro}

This project was an attempt to explore the capabilities and limitations of large language models (LLMs) in research-level mathematics.
It all started at the \href{https://combinatorial-cowork.space/editions/2026/}{Combinatorial Coworkspace} in March 2026, where the three event organizers decided to invite research mathematicians to contribute and discuss benchmark questions using the \href{https://math.sciencebench.ai}{ScienceBench} platform:
\vspace*{-3pt}
\begin{center}
    \url{https://math.sciencebench.ai/benchmarks/2026-09-28}\,.
\end{center}

This document reports on this compilation and in particular on the \emph{Benchmarks in Leipzig} workshop, which took place May 11--13, 2026, at the \href{https://www.mis.mpg.de/}{Max Planck Institute for Mathematics in the Sciences} in Leipzig, Germany:
\vspace*{-3pt}
\begin{center}
    \url{https://www.mis.mpg.de/events/series/benchmarks-in-leipzig}\,.
\end{center}

\begin{findingsbox}
\begin{itemize}
    \item Writing benchmark questions in the style of exercises based on publicly accessible research has reached its limits when it comes to the best-performing available models, see \Cref{tab:kof3-gptpro,tab:stage4-evals}.
    For all other models, this approach remains informative.
    Between Stages~1 and~4, the number of questions solved in a single run by the best model rose from~44 to~91, the latter with web search and code execution enabled.
    \item The considered LLMs differ drastically in performance, see \Cref{tab:single-run-evals,,tab:kof20-stage2}.
    \item The performance of each model in \Cref{tab:kof20-stage2} varies greatly between runs.
    \item A model may perform very differently with and without code execution.
    Without a code interpreter, \sol exhausted its token budget in repeated web searches and solved only~68 of the~100 questions, while it solved~90 with a code interpreter enabled, see \Cref{sec:method:rerun}.
    \item Many of the researchers involved used the ScienceBench chat, which suggests that research mathematicians are interested in using AI chat tools when given access to reasoning models in a trusted environment.
    \item The LLM review of submissions and the comparison with model-generated solutions led to corrections in several questions, and to three questions being removed from the benchmark.    
    This indicates the potential benefit of AI support in the peer review of mathematical content.
    \item A single strong attempt does not replace many weak ones: two of the four questions that no model solved in Stage~4 had already been solved in Stage~2, by weaker models that were given~20 attempts instead of one.
\end{itemize}
\end{findingsbox}

During the question contribution period, April 1--May 15, 2026,~\textbf{49 research mathematicians} contributed benchmark questions, and~35 of them participated in person in the Benchmarks in Leipzig workshop.
During the workshop, the participants combined research exploration via the \textbf{ScienceBench chat interface} with structured benchmark submission.
The per-user chat has the same project-active models that run on the benchmark, see \Cref{tab:single-run-evals} and also \Cref{sec:platform}.
During the project phase,~39 of the contributors used it, \textbf{accumulating~1,067 turns across~483 chats}.
The researchers contributed a total of~\textbf{100 research-level questions} with unique, well-defined answers.
Participants were encouraged to create questions that were based on existing and publicly accessible research but went beyond what could be presented in a research paper.
Unlike in the \emph{First Proof} project by~\cite{firstproof}, the goal was not to provide solutions to currently unpublished open problems.

The~100 benchmark questions are presented in Appendix~\ref{app:problems} and belong to the following (overlapping) research areas:

\begin{table}[h!]
\centering
\begin{tabular}{@{}lc@{\hspace{1em}}|lc@{\hspace{1em}}|lc@{}}
\toprule
Research area &  & Research area &  & Research area &  \\
\midrule
  Algebraic Geometry & 36 & Commutative Algebra & 5 & Topology & 3 \\
  Algebraic Combinatorics & 22 & Algebraic Statistics & 4 & Complex Analysis & 2 \\
  Matroid Theory & 18 & Homological Algebra & 4 & Euclidean Geometry & 2 \\
  Representation Theory & 15 & Polytope Theory & 4 & Metric Geometry & 1 \\
  Enumerative Combinatorics & 15 & Tropical Geometry & 4 & Probability Theory & 1 \\
  Combinatorics & 14 & Analysis & 3 & Real Algebraic Geometry & 1 \\
  Discrete Geometry & 10 & Arithmetic Geometry & 3 & Theoretical Computer Science & 1 \\
  Algebra & 7 & Knot Theory & 3 & & \\
  Graph Theory & 5 & Number Theory & 3 & & \\
\bottomrule
\end{tabular}
\vspace*{3pt}
\caption{Research areas of the questions (overlapping), generated using \claude and subsequently verified by the authors.}
\vspace*{-10pt}
\end{table}

\paragraph{Evaluation}\

\smallskip
It turned out that most of our questions could be solved by at least one state-of-the-art model, though the differences between models, and between repeated runs of the same model, were substantial.
The evaluation process for each question consisted of four stages, summarized in \Cref{tab:stages}:

\begin{table}[h!]
\centering
\begin{tabular}{@{}cccccc@{}}
\toprule
Stage & Date & Models & Attempts & Best & Unsolved \\
\midrule
 1 & May 13 & 10 &  1 & 44 & 41 \\
 2 & May 22 &  3 & 20 & 75 & 16 \\
 3 & Jun 2  &  2 &  3 & 88 &  2 \\
 4 & Sep 28 &  9 &  1 & 91 &  1 \\
\bottomrule
\end{tabular}
\vspace*{3pt}
\caption{The four evaluation stages. In Stage~4, two of the seven models were also run in their Pro mode, giving nine configurations.\label{tab:stages}}
\vspace*{-10pt}
\end{table}

Here \emph{Attempts} is how often each model was asked each question, \emph{Best} the number of questions solved by the strongest model of a stage, and \emph{Unsolved} the number of questions that no model had solved in that stage or in any earlier one.
The four stages cannot be compared directly, because they differ in the models used, in the number of attempts per question, and in whether the models could access the internet and execute code.

\textbf{Stage~1:}
During the submission process, the following five \emph{project-active models} (left) and their predecessors (right) each attempted the question in a single run\footnote{If a model timed out or produced an error, we reran the model up to three times. Further, we originally had \texttt{Grok 4.20} active instead of \grok; we reran all questions through \grok afterward for an up-to-date comparison.}.
The following table shows the models' results:

\begin{table}[h!]
\centering
\begin{minipage}[t]{0.48\linewidth}
\vspace*{0pt}\centering
\begin{tabular}{@{}lrc@{}}
\toprule
Model & Release date & Solved \\
\midrule
  \gpt & 2026-04-23 & 44 \\
  \gemini & 2026-03-17 & 15 \\
  \claude & 2026-04-15 & 13 \\
  \deepseek & 2026-04-24 & 10 \\
  \grok & 2026-05-06 & 6 \\
\bottomrule
\end{tabular}
\end{minipage}\hfill
\begin{minipage}[t]{0.48\linewidth}
\vspace*{0pt}\centering
\begin{tabular}{@{}lrc@{}}
\toprule
Model & Release date & Solved \\
\midrule
  \texttt{GPT-5.4} & 2026-03-17 & 28 \\
  \texttt{Gemini 3 Pro} & 2025-11-18 & 14 \\
  \texttt{Claude Opus 4.6} & 2026-02-04 & 14 \\
  \texttt{DeepSeek-V3.2} & 2025-12-01 & 8 \\
  \texttt{Grok 4.20} & 2026-02-20 & 5 \\
\bottomrule
\end{tabular}
\end{minipage}
\vspace*{3pt}
\caption{Release dates and results of the project-active models (left) and their predecessors (right).\label{tab:single-run-evals}}
\vspace*{-10pt}
\end{table}

In Stage~1,~\textbf{41 questions remained unsolved} by all models.

\textbf{Stage~2:}
After the workshop on May 11--13, \textbf{Surge AI} queried on our behalf each of three models~20 times on each of the~100 questions.
\Cref{tab:surge-leaderboard} summarizes the results.
We say that a model \emph{answered} a question during a run if it produced any answer at all.
Because each question was run multiple times per model, we distinguish two metrics: a \emph{correct run} is a single model attempt that produced the correct answer, while a \emph{solved question} is one for which at least one of the~20 attempts was a correct run.

\begin{table}[h!]
\centering
\begin{tabular}{@{}lcrrc@{}}
\toprule
Model & Runs & Answered & Correct runs & Solved questions \\
\midrule
  \gpt    & 2,000 & 1,389 (69\,\%) & 1,043 (52\,\%) & 75 (75\,\%) \\
  \gemini & 2,000 &  832 (42\,\%) &  406 (20\,\%) & 40 (40\,\%) \\
  \claude & 2,000 & 1,306 (65\,\%) &  294 (15\,\%) & 44 (44\,\%) \\
\bottomrule
\end{tabular}
\vspace*{3pt}
\caption{Results of three models based on 20 external runs per question.\label{tab:surge-leaderboard}
}
\vspace*{-10pt}
\end{table}

Observe that \claude had rather few correct runs compared with its number of solved questions.
The reason is that it solved many questions correctly only a few times.
For instance, there were 19 questions that it solved in no more than 3 of its 20 attempts.
We refer to \Cref{sec:crossmodelevals} for additional performance statistics.

In Stage~2, \textbf{20 questions remained unsolved} by all models, and in Stages~1--2, \textbf{16 questions remained unsolved}.

\textbf{Stage~3:}
All~100 questions were then attempted three times each by \texttt{GPT-5.5 Pro} with Extended Thinking enabled, which we ran through the \texttt{ChatGPT} website, where web search and code execution are available.
In addition, \textbf{Surge AI} queried \texttt{Gemini 3.1 Pro Deep Think} three times on each question.
The following table summarizes the results:

\begin{table}[h!]
\centering
\begin{tabular}{@{}lcrrc@{}}
\toprule
Model & Runs & Answered & Correct runs & Solved questions \\
\midrule
  \texttt{GPT-5.5 Pro} {\small (Extended Thinking)}    & 300 & 300 (100\,\%) & 237 (79\,\%) & 88 (88\,\%) \\
  \texttt{Gemini 3.1 Pro Deep Think} & 300 & 250 (83\,\%) & 141 (47\,\%) & 56 (56\,\%) \\
\bottomrule
\end{tabular}
\vspace*{3pt}
\caption{Results of the two heavy-thinking models across 3 runs each.\label{tab:stage3-leaderboard}
}
\vspace*{-10pt}
\end{table}

Note that the two rows cannot be compared directly:
when \texttt{GPT-5.5 Pro} failed to produce an answer to a question at all, we reran it multiple times on the \texttt{ChatGPT} website.
On the other hand, we did not rerun \texttt{Gemini 3.1 Pro Deep Think} on the questions it failed to answer, which is why it produced no answer on 14, 17, and 19 of the 100 questions across its three runs, respectively.
We refer to \Cref{sec:crossmodelevals} for additional performance statistics.

In Stage~3, \textbf{8 questions remained unsolved} by all models, and in Stages~1--3, \textbf{2 questions remained unsolved}.

\textbf{Stage~4:}
Four months after the workshop, in September 2026, all~100 questions were attempted once more by the models that had become available in the meantime, seven models in nine configurations (two of the models were also run in their Pro mode)\footnote{We did not include \texttt{Grok 4.7}, released on September 21, 2026, because we could not bound its usage: it answers by running its own tool loop, and a single test run consumed more than 9 million input tokens.}, each in a single run\footnote{As in Stage~1, we reran a model that timed out or produced an error. \gemflash was run without code execution, because its code interpreter regularly exceeds the provider's limit on tool calls, and it was allowed a second attempt. It still produced no answer on~41 of the~100 questions, all of which are counted as unsolved.}.
Unlike in Stage~1, web search and code execution were enabled, see \Cref{sec:method:rerun}.
The following table shows the models' results:

\begin{table}[h!]
\centering
\begin{tabular}{@{}lrc@{}}
\toprule
Model & Release date & Solved \\
\midrule
  \astra    & 2026-09-03 & 91 \\
  \solpro   & 2026-07-09 & 91 \\
  \sol      & 2026-07-09 & 90 \\
  \astrapro & 2026-09-03 & 88 \\
  \gptsixsol   & 2026-09-22 & 86 \\
  \opusfivefive & 2026-09-22 & 85 \\
  \fable    & 2026-09-01 & 79 \\
  \opusfive & 2026-07-24 & 75 \\
  \gemflash & 2026-09-02 & 41 \\
\bottomrule
\end{tabular}
\vspace*{3pt}
\caption{Release dates and results of the models in the September 2026 re-run.\label{tab:stage4-evals}}
\vspace*{-10pt}
\end{table}

In Stage~4, \textbf{4 questions remained unsolved} by all models, and in Stages~1--4, \textbf{1 question remained unsolved}.

\medskip

This benchmark sits alongside other recent benchmarks for research-level mathematics, such as FrontierMath~\cite{frontiermath}, First Proof~\cite{firstproof}, IMProofBench~\cite{improofbench}, Soohak~\cite{soohak}, RealMath~\cite{realmath}, and Riemann-Bench~\cite{riemannbench}.

\section{Building the Leipzig Benchmark}
\label{sec:method}

The contributors jointly compiled the Leipzig Benchmark between April 1 and May 15, 2026, with the bulk of the questions being submitted during the Benchmarks in Leipzig workshop at the MPI-MiS Leipzig, see \Cref{sec:event}.
The three event organizers led the invitation process for the contributions.
We aimed to have 35 on-site contributors, and then decided to additionally invite~14 external contributors.

\subsection{Submission guidelines}
\label{sec:method:guideline}

Contributors were invited to submit their questions via the ScienceBench platform, see \Cref{sec:platform}.
They were given the following guidelines:
\begin{itemize}
  \item The question must have a unique, well-defined, unguessable answer.
  \begin{itemize}
      \item Possible answers include, for example, a number, a polynomial, a counting function with parameters, or an expression or formula.
      \item The answer to the question should not be a proof, ``necessary'' or ``sufficient'' conditions, or any similar derivation.
      \item The ``unguessable'' condition roughly means that a reader allowed to guess the answer multiple times without working out any details should not have a realistic chance of guessing correctly.
  \end{itemize}
  \item Make the question difficult enough so that only a few of the five project-active models (see \Cref{tab:single-run-evals}) can solve it,
  if any.
  \begin{itemize}
      \item During the submission workflow (see \Cref{sec:method:workflow}), only questions answered correctly by at most three of these models were accepted for final submission.
  \end{itemize}
  \item The complexity of the question should not arise solely from a computation.
  \begin{itemize}
      \item Computer algebra systems may be used in solving the question, but such a computation should not be at its core.
  \end{itemize}
  \item The models should not fail just because
    (a) they base their answers on a trivially wrong assumption (for example, misinterpreting a variable called~$\pi$ as the number~$\pi$), or
    (b) they have not yet learned the source material.
  This means, in particular, that the solution to the question should not use results from work that is not publicly accessible.
\end{itemize}

\paragraph{Comments}
\begin{itemize}
  \item A guiding idea was to create questions that are based on existing research but go beyond what can be presented in a research paper (for instance, a lengthy example that does not fit into an article).
  \item The goal was not to pose open problems, nor questions whose solutions rely on proofs that have not yet been published.
  \item The rule not to rely on unpublished work turned out to be a notable constraint.
  Many contributors would have preferred to submit questions related to their current research.
  \item Unguessability was sometimes hard to achieve, as several questions turned out to have only a very small number of possible answers, which could be guessed rather easily.
\end{itemize}

\subsection{Submission workflow}
\label{sec:method:workflow}

The submission workflow on the ScienceBench platform consisted of three stages:
\begin{itemize}[labelwidth=1.5cm,leftmargin=2cm]
  \item[Stage~W1] The contributor entered a question together with its answer.
  \item[Stage~W2] The five project-active models, see \Cref{tab:single-run-evals}, attempted the question.
  Their answers were compared with the original answer by an LLM judge.
  Only if at most three of the five models answered correctly could the contributor proceed.
  \item[Stage~W3] The contributor entered a worked-out solution to the question that resulted in the already-provided answer, as well as relevant references.
\end{itemize}

\paragraph{Comments}
\begin{itemize}
  \item The solution and the references played no role in the evaluation: they were shown neither to the models nor to the other contributors.
  Their sole purpose was to make the contributors think through their question and their derivation of the answer carefully.
  The answer, too, was withheld from the five models attempting the question; it was seen only by the LLM judges.
  \item In several cases, writing the solution down revealed a mistake in the original construction and led to a corrected answer.
  \item In practice, many contributors revised their questions, and sometimes their answers, in light of the models' replies, even when the question had already passed the threshold of at most three models answering it correctly.
\end{itemize}

\subsection{Post-collection audit}
\label{sec:method:audit}

The benchmark questions went through three audit stages:

\begin{itemize}[labelwidth=1.5cm,leftmargin=2cm]
  \item[Stage~A1] All contributors were able to see all questions and the corresponding model answers and solutions.
  In several cases, this triggered discussions that resulted in fixing minor or even major errors in the originally submitted question/answer.
  \item[Stage~A2] After the workshop ended, each question went through an LLM review.
  Each submission was reviewed for typos in the mathematics, inconsistencies, missing content, and possible reasoning errors.
  \item[Stage~A3] The questions that remained completely unsolved by all models after the evaluation in Stage~3 went through another review by the original authors.
\end{itemize}

\paragraph{Comments}
\begin{itemize}
  \item The AI-assisted review of the Leipzig benchmark at the beginning of Stage~A2 flagged~16 potential mathematical issues in the submissions and triggered personal email exchanges with their authors, after which the submissions were updated accordingly.
  This led to~12 (mostly minor) revisions of question statements, and~3 questions were subsequently removed from the benchmark.
  \item This feedback loop continued during the evaluation rounds in Stages~2 and~3 and resulted in~4 further answers being corrected by the authors,~3 of them to the answer that \texttt{GPT-5.5 Pro} had given and~1 to something different from the model's answer.
  \emph{The numbers presented in the introduction already reflect these updates.}
  \item Such modifications to benchmark problem sets after an LLM review are not unusual.
  Epoch AI, the company behind the \emph{FrontierMath} project, which comprises mathematics problems of all difficulties, announced on May 11, 2026, that an AI-assisted review of its problems in Tiers~1--4 had flagged fatal errors in about a third of the questions~\cite{epochai_frontiermath_audit_2026}.
  The corrected set appeared a month later: 135 of the 350 problems had been revised and a further 12 removed, so that errors had been found in 42\,\% of them~\cite{epochai_frontiermath_v2_2026}.

\end{itemize}

\subsection{Single-attempt model runs}
\label{sec:method:models}

In Stage~W2 of the submission workflow, the following model runs were involved, see \Cref{tab:single-run-evals} for evaluation:

\begin{itemize}[labelwidth=1.5cm,leftmargin=2cm]
  \item[Check] Before a question was sent to the models to solve, it was checked by \gpt and \gemini against the submission guidelines.
  \item[Runs] After passing the check, the question was attempted by the five project-active models.
  The strongest available version of each model was queried via the API; the exact configurations are listed in \Cref{tab:api-settings}.
  \item[Answer] After the models finished, their answers were compared with the original answer using \gpt and \gemini as LLM judges.
\end{itemize}

\paragraph{Model settings}\ 

\smallskip
The API configurations used in Stage~W2 to query the five project-active models are presented in the following table:

\begin{table}[h!]
\centering
\begin{tabular}{@{}llc@{}}
\toprule
Model & Effort & Maximum Token Usage \\
\midrule
  \gpt       & none / low / medium / high / \textbf{xhigh} & 128k / 128k \\
  \gemini    & minimal / low / \textbf{medium} / high & 64k / 64k \\
  \claude    & low / medium / \textbf{high} / xhigh / max & 128k / 128k \\
  \deepseek  & \textbf{high} / max & 128k / 384k \\
  \grok      & none / low / medium / \textbf{high} & [no parameter] \\
\bottomrule
\end{tabular}
\vspace*{3pt}
\caption{API settings used for the project-active models.}
\label{tab:api-settings}
\end{table}
\vspace*{-5pt}

The effort used in our API calls is shown in bold among the available effort options for each model.
The third column lists the maximum allowed token usage (left value) next to the models' actual maximum (right value).

\paragraph{Comments}
\begin{itemize}
  \item The complexity of the questions caused many model runs to end in an error or a timeout.
  We then reran these models up to three times on such questions; after three failed attempts, we recorded the result as incorrect.
  \item Web search, code execution, and computation tools were all disabled for the solve runs in Stage~W2.
  Earlier runs with code execution and computational tools enabled had shown that the models resorted to brute-force computations, which resulted in timeouts and errors, and the models in general performed worse than they did when code execution was disabled.
  One example from the ``chain of thought'' of \gpt: \emph{I'll translate your question into the [...] framework, then compute the answer exactly via [...] rather than doing simulation.}
  \item We chose to run \gemini, \claude, and \deepseek\ with submaximal effort because otherwise these models timed out too frequently.
  Moreover, the maximal token usage for \deepseek\ is set below its model maximum to match the cap used for the other models.
\end{itemize}

\subsection{Multi-run evaluations}
\label{sec:method:external}

Stages~2 and~3 consisted of the following multi-run evaluations:
\begin{itemize}[labelwidth=1.5cm,leftmargin=2cm]
  \item[Stage~2] All~100 questions were sent to \textbf{Surge AI}, which queried each of the three models in \Cref{tab:surge-leaderboard}~20 times per question.
  \item[Stage~3] All~100 questions were then attempted three times each by \texttt{GPT-5.5 Pro} with Extended Thinking enabled, which we ran through the \texttt{ChatGPT} website, where web search and code execution are available.
  Additionally, \textbf{Surge AI} queried \texttt{Gemini 3.1 Pro Deep Think} three times on each question.
\end{itemize}
The overall results of Stages~2 and~3 are reported in \Cref{tab:surge-leaderboard,,tab:stage3-leaderboard} in the \hyperref[sec:intro]{Introduction}, respectively.

\subsection{September 2026 re-run}
\label{sec:method:rerun}

Stage~4 repeated the format of the single-run setup of Stage~1 with the new models that had become available by September 2026, see \Cref{tab:stage4-evals}:

\begin{itemize}[labelwidth=1.5cm,leftmargin=2cm]
  \item[Runs] Each of the nine configurations attempted each of the~100 questions once, queried through the models' APIs.
  \item[Tools] Web search and code execution were enabled, with the exception of \gemflash, for which code execution had to be disabled.
  \item[Effort] The five configurations for the GPT models ran at maximal effort with a cap of 128k tokens, and \opusfive, \fable, and \opusfivefive ran at high effort with an advisory budget of 100k tokens, as maximal effort produced the same answers using twice the tokens.
\end{itemize}

\paragraph{Comments}
\begin{itemize}
  \item Enabling code execution, which we had switched off for the solve runs in Stage~1, had a substantial effect on the results.
  The same model, \sol, solved only 68 of the~100 questions at xhigh effort, with web search but without a code interpreter, whereas it answered every question and solved~90 when code execution was enabled and it ran at maximal effort.
  Without a code interpreter, it frequently exhausted its token budget in repeated web searches and produced no answer at all on~16 of the~100 questions.
  The newer models of Stage~4 balance code execution and web search with reasoning much better than the models of Stage~1 did.
  \item Many of the questions lean on results from recent papers by their contributors; these papers have since been, or will soon be, published and indexed, so future training corpora will contain that material.
  Moreover, the questions (though not their answers) have been public since June 2026, and all models of Stage~4 were released after that.
  We therefore expect future models to do better on this collection than the present ones do, for reasons that have nothing to do with mathematical ability, and this may already be visible in Stage~4.
  This is a structural property of any benchmark built at the research frontier, and it is the reason we consider the dated evaluations of \Cref{tab:stages}, rather than a single number, to be the result.
\end{itemize}

\subsection{Cross-run and cross-model statistics}\
\label{sec:crossmodelevals}

Beyond the aggregated statistics in \Cref{tab:surge-leaderboard,tab:stage3-leaderboard,tab:stage4-evals}, we provide the following refined counts.
We start with the distributions of correct solutions produced during the~20 runs per model in Stage~2.

\begin{table}[h!]
\centering
\begin{tabular}{lcccccc}
\toprule
& \multicolumn{6}{c}{Number of correct runs (out of 20)} \\
\cmidrule(lr){2-7}
Model & 0 & 1 to 4 & 5 to 9 & 10 to 14 & 15 to 19 & 20 \\
\midrule
\gpt    & 25 & 12 &  9 &  7 & 20 & 27 \\
\gemini & 60 & 12 &  8 &  5 & 10 &  5 \\
\claude & 56 & 20 &  8 & 10 &  5 &  1 \\
\bottomrule
\end{tabular}
\vspace*{5pt}
\caption{Per-question consistency distribution across the 20 runs in Stage~2.}
\label{tab:kof20-stage2}
\vspace*{-10pt}
\end{table}

Table~\ref{tab:kof3-gptpro} shows the analogous distributions of correct solutions in the~3 runs in Stage~3.
\begin{table}[h!]
\centering
\begin{tabular}{l >{\centering\arraybackslash}p{30pt} >{\centering\arraybackslash}p{30pt} >{\centering\arraybackslash}p{30pt} >{\centering\arraybackslash}p{30pt}}\toprule
& \multicolumn{4}{c}{Number of correct runs (out of 3)} \\
\cmidrule(lr){2-5}
Model & 0 & 1 & 2 & 3 \\
\midrule
\texttt{GPT-5.5 Pro} {\small (Extended Thinking)} & 12 & 8 & 11 & 69 \\
\texttt{Gemini 3.1 Pro Deep Think} & 44 & 10 & 7 & 39 \\
\bottomrule
\end{tabular}
\vspace*{5pt}
\caption{Per-question consistency distribution across the 3 runs in Stage~3.}
\label{tab:kof3-gptpro}
\vspace*{-10pt}
\end{table}
We continue with a cross-comparison of the two heavy-thinking models, shown in Table~\ref{tab:pro-vs-dt}.
\begin{table}[h!]
\centering
\small
\begin{tabular}{llcc}
\toprule
& & \multicolumn{2}{c}{\texttt{Gemini 3.1 Pro Deep Think}} \\
\cmidrule(lr){3-4}
& & Solved & Unsolved \\
\midrule
\multirow{2}{*}{\texttt{GPT-5.5 Pro} {\small (Extended Thinking)}} & Solved   & 52 & 36 \\
                                         & Unsolved &  4 &  8 \\
\bottomrule
\end{tabular}
\vspace*{5pt}
\caption{Contingency table cross-comparing the two models in Stage~3.}
\label{tab:pro-vs-dt}
\vspace*{-10pt}
\end{table}

We conclude with a cross-comparison of all nine configurations of Stage~4.
The diagonal entries are the numbers of questions solved by the respective configuration, and the off-diagonal entries count the questions solved by both configurations.
Of the~100 questions,~31 were solved by every one of the nine and~96 by at least one of them, and no question was solved by only a single configuration:

\begin{table}[h!]
\centering
\small
\setlength{\tabcolsep}{4pt}
\begin{tabular}{@{}lccccccccc@{}}
\toprule
 & \shortstack{\texttt{GPT-6}\\\texttt{Astra}} & \shortstack{\texttt{GPT-5.6}\\\texttt{Sol Pro}} & \shortstack{\texttt{GPT-5.6}\\\texttt{Sol}} & \shortstack{\texttt{GPT-6}\\\texttt{Astra Pro}} & \shortstack{\texttt{GPT-6}\\\texttt{Sol}} & \shortstack{\texttt{Claude}\\\texttt{Opus 5.5}} & \shortstack{\texttt{Claude}\\\texttt{Fable 5.1}} & \shortstack{\texttt{Claude}\\\texttt{Opus 5}} & \shortstack{\texttt{Gemini}\\\texttt{3.8 Flash}} \\
\midrule
\astra         & 91 & 87 & 86 & 85 & 86 & 81 & 74 & 72 & 40 \\
\solpro        & 87 & 91 & 89 & 83 & 83 & 83 & 78 & 74 & 40 \\
\sol           & 86 & 89 & 90 & 82 & 82 & 83 & 77 & 72 & 40 \\
\astrapro      & 85 & 83 & 82 & 88 & 81 & 78 & 71 & 70 & 39 \\
\gptsixsol     & 86 & 83 & 82 & 81 & 86 & 76 & 71 & 68 & 37 \\
\opusfivefive  & 81 & 83 & 83 & 78 & 76 & 85 & 77 & 73 & 40 \\
\fable         & 74 & 78 & 77 & 71 & 71 & 77 & 79 & 70 & 38 \\
\opusfive      & 72 & 74 & 72 & 70 & 68 & 73 & 70 & 75 & 37 \\
\gemflash      & 40 & 40 & 40 & 39 & 37 & 40 & 38 & 37 & 41 \\
\bottomrule
\end{tabular}
\vspace*{5pt}
\caption{Cross-comparison of all nine configurations in Stage~4.}
\label{tab:stage4-overlap}
\vspace*{-10pt}
\end{table}
\section{The ScienceBench platform}
\label{sec:platform}

The benchmark is hosted on the ScienceBench platform at \url{https://math.sciencebench.ai}.
The platform consists of publicly available benchmarks and several components that are accessible to registered users, which include:
\begin{itemize}
    \item a three-stage submission system described in \Cref{sec:method:workflow},
    \item a browsing interface that lets every registered user view other contributors' questions, comment on them, attempt the questions, and inspect each model's full answer side by side, and
    \item a chat system that allows users to chat with all project-active models in parallel.
\end{itemize}

\section{The Benchmarks in Leipzig workshop}
\label{sec:event}

The bulk of the benchmark was contributed during the three-day Benchmarks in Leipzig workshop held at the Max Planck Institute for Mathematics in the Sciences on May 11--13, 2026.
Each day consisted of three phases:
\begin{enumerate}
    \item joint discussions or presentations on AI methods and tools in mathematics research;
    \item discussion and work in small groups;
    \item wrap-up: a showcase of some submissions, a summary of progress, and a discussion of the difficulties the participants had encountered in crafting their questions.
\end{enumerate}
The groups were shuffled each day to encourage interaction and collaboration among researchers from different areas of mathematics.
This also allowed participants to exchange experiences gained during earlier sessions.

\begin{itemize}[labelwidth=1.5cm,leftmargin=2cm]
\item[Day~1]
The workshop opened with a short introduction by Christian Stump to the ScienceBench platform (\Cref{sec:platform}), followed by an overview of the submission guidelines (\Cref{sec:method:guideline}).
The participants then worked in preassigned groups of two to four, each group composed to mix scientific backgrounds and degrees of prior familiarity with ScienceBench and AI research tools.

\item[Day~2]
The second day began with a presentation by Nupur Jain on the use of machine learning methods in mathematics research.
This was followed by a discussion led by Bernd Sturmfels on the extent to which AI tools are, can be, and should be used in mathematics research, and on how they are shaping the field.
Throughout the day, all participants worked on and discussed their submissions, sometimes using the ScienceBench chat system to do so.
Toward the end of the day, a few submissions were presented by their authors and discussed jointly.

\item[Day~3]
After a brief discussion of the status of the submissions, the remainder of the workshop was devoted to finalizing drafts, solving and reviewing other participants' submissions, and further exploring the capabilities of the ScienceBench chat system.
Some of the newly submitted questions were presented, and the models' solutions were compared with the contributors' intended approaches.
By the end of the workshop, and after the post-collection audit (\Cref{sec:method:audit}), the Leipzig Benchmark consisted of~100 submissions meeting the submission guidelines, contributed by~35 on-site and~14 off-site contributors.
\end{itemize}

\appendix
\section{The Leipzig Benchmark questions}
\label{app:problems}

The 100 Leipzig Benchmark questions are listed below.
Each box title carries a label in italics that records the stage at which the question was first answered correctly:

\begin{center}
\begin{tabular}{@{}l@{\hspace{2em}}l@{\hspace{1em}}r@{}}
\toprule
Label & Meaning & Count \\
\midrule
\emph{solved in Stage 1} & correct answer produced by one of the 5 project-active models & 59 \\
\emph{solved in Stage 2} & first solved by one of Surge AI's multi-runs & 25 \\
\emph{solved in Stage 3} & first solved by one of the multi-runs of the two heavy-thinking models & 14 \\
\emph{remains unsolved} & unsolved by all of the above & 2 \\
\bottomrule
\end{tabular}
\end{center}

These labels record the outcome of Stages~1 to~3.
In Stage~4, the number of unsolved questions changed from~2 to~1.

\input{problems.tex}

\clearpage
\section*{Affiliations}
\label{app:affiliations}

\begin{enumerate}\itemsep0pt
  \item Ruhr University Bochum, Bochum, Germany
  \item Max Planck Institute for Mathematics in the Sciences, Leipzig, Germany
  \item TU Berlin, Berlin, Germany
  \item University of Florida, Gainesville, FL, USA
  \item Ecole Normale Superieure (ENS), PSL University, Paris, France
  \item University of California, Davis, USA
  \item University of Barcelona, Barcelona, Spain
  \item University of California, Berkeley, USA
  \item ETH Zurich, Zurich, Switzerland
  \item Universidad de Talca, Talca, Chile
  \item TU Dresden, Dresden, Germany
  \item University of Bonn, Bonn, Germany
  \item Universit\'e du Qu\'ebec \`a Montr\'eal (UQAM), Montreal, Canada
  \item University of Southern California, Los Angeles, USA
  \item University of Luxembourg, Luxembourg
  \item Amherst College, Amherst, Massachusetts, USA
  \item Bielefeld University, Bielefeld, Germany
  \item Technische Universit\"at Braunschweig, Braunschweig, Germany
  \item University of Leipzig, Leipzig, Germany
  \item University of Copenhagen, Copenhagen, Denmark
  \item KTH Royal Institute of Technology, Stockholm, Sweden
  \item Universit\'e de Montr\'eal, Montr\'eal, Canada
  \item Georg-August-Universit\"at G\"ottingen, G\"ottingen, Germany
  \item Goethe University Frankfurt, Frankfurt am Main, Germany
  \item University of Texas at Dallas, Richardson, Texas, USA
\end{enumerate}

\end{document}

%% file: problems.tex

\begin{problembox}[Question 001\hfill\emph{solved in Stage 1}]
Let $G$ be a matrix in $F_3^{4\times n}$ with rows such that for every vector list of 287 vectors $v_1,\ldots,v_{287}$ (there can be repetition) in ${F_3}^4$ I can find $287$ disjoint subsets $X_1,\dots,X_{287}$ of columns of $G$ such that $v_i$ belongs to the span of the columns of $X_i$. Let $C$ be a linear code generated by a matrix with the property described above that minimise $n$. What is the weight enumerator of a code generated by such a $G$?
\end{problembox}

\begin{problembox}[Question 002\hfill\emph{solved in Stage 1}]
Consider the following homogeneous ideal in $\mathbb{C}[z_1,\dots,z_5]$:
\begin{multline*}
\bigl\langle -z_{2}+z_{4}-z_{5},\\
\quad z_{1}^{5}z_{2}^{5}+5z_{1}^{4}z_{2}^{5}z_{3}+10z_{1}^{3}z_{2}^{5}z_{3}^{2}+10z_{1}^{2}z_{2}^{5}z_{3}^{3}+5z_{1}z_{2}^{5}z_{3}^{4}-z_{2}^{6}z_{3}^{4}+z_{2}^{5}z_{3}^{5}\\
+5z_{1}^{5}z_{2}^{4}z_{4}+25z_{1}^{4}z_{2}^{4}z_{3}z_{4}+50z_{1}^{3}z_{2}^{4}z_{3}^{2}z_{4}+50z_{1}^{2}z_{2}^{4}z_{3}^{3}z_{4}+25z_{1}z_{2}^{4}z_{3}^{4}z_{4}+5z_{2}^{4}z_{3}^{5}z_{4}\\
+10z_{1}^{5}z_{2}^{3}z_{4}^{2}+50z_{1}^{4}z_{2}^{3}z_{3}z_{4}^{2}+100z_{1}^{3}z_{2}^{3}z_{3}^{2}z_{4}^{2}+100z_{1}^{2}z_{2}^{3}z_{3}^{3}z_{4}^{2}+50z_{1}z_{2}^{3}z_{3}^{4}z_{4}^{2}+10z_{2}^{3}z_{3}^{5}z_{4}^{2}\\
+10z_{1}^{5}z_{2}^{2}z_{4}^{3}+50z_{1}^{4}z_{2}^{2}z_{3}z_{4}^{3}+100z_{1}^{3}z_{2}^{2}z_{3}^{2}z_{4}^{3}+100z_{1}^{2}z_{2}^{2}z_{3}^{3}z_{4}^{3}+50z_{1}z_{2}^{2}z_{3}^{4}z_{4}^{3}+10z_{2}^{2}z_{3}^{5}z_{4}^{3}\\
+5z_{1}^{5}z_{2}z_{4}^{4}+25z_{1}^{4}z_{2}z_{3}z_{4}^{4}+50z_{1}^{3}z_{2}z_{3}^{2}z_{4}^{4}+50z_{1}^{2}z_{2}z_{3}^{3}z_{4}^{4}+25z_{1}z_{2}z_{3}^{4}z_{4}^{4}+5z_{2}z_{3}^{5}z_{4}^{4}\\
+z_{1}^{5}z_{4}^{5}+5z_{1}^{4}z_{3}z_{4}^{5}+10z_{1}^{3}z_{3}^{2}z_{4}^{5}+10z_{1}^{2}z_{3}^{3}z_{4}^{5}+5z_{1}z_{3}^{4}z_{4}^{5}+z_{3}^{5}z_{4}^{5}\\
-6z_{2}^{5}z_{3}^{4}z_{5}-15z_{2}^{4}z_{3}^{4}z_{5}^{2}-20z_{2}^{3}z_{3}^{4}z_{5}^{3}-15z_{2}^{2}z_{3}^{4}z_{5}^{4}-6z_{2}z_{3}^{4}z_{5}^{5}-z_{3}^{4}z_{5}^{6}\bigr\rangle.
\end{multline*}
Can you tell me the Picard group of the projective variety cut out by this ideal?
\end{problembox}

\begin{problembox}[Question 003\hfill\emph{solved in Stage 1}]
What is the last fall degree of the family of polynomials $\{x^{12456789}+y,y^{12456789}+x,x^{503}y^{503}\}$?
\end{problembox}

\begin{problembox}[Question 004\hfill\emph{solved in Stage 1}]
Consider the ideal $I \subset \mathbb{Q}[x,y]$ generated by $5$ polynomials:
\begin{align*}
I = \bigl\langle\;
& -1 - ax + x^3 - by + y^2,\\
& -c - dx^2 - y + xy - ey^2,\\
& -fx + x^2 - gy + x^2 y - y^2,\\
& -h + x - kx^2 - my - ny^2 + xy^2,\\
& -p - qx - rx^2 + (57y)/13 + y^3
\;\bigr\rangle,
\end{align*}
where $(a, b, c, d, e, f, g, h, k, m, n, p, q, r)$ are rational numbers. There is one set of values for $(a, b, c, d, e, f, g, h, k, m, n, p, q, r)$ for which the quotient $\mathbb{Q}[x,y]/I$ has the structure of a $5$-dimensional $\mathbb{Q}$-vector space with basis given by $\{[1],[x],[x^2],[y],[y^2]\}$. Using those values, find $m+n$.
\end{problembox}

\begin{problembox}[Question 005\hfill\emph{solved in Stage 1}]
Let $A$ be the connected cyclic Nakayama algebra with Kupisch series $[2,3]$ and let $B$ be the Auslander algebra of $A$. Let $C$ be the algebra $B/I$ where $I$ is the two-sided ideal given by the socle of the direct sum of all indecomposable projective-injective $B$-modules (as a submodule of the regular module). What is the sum of all delooping levels of the indecomposable $C$-modules?
\end{problembox}

\begin{problembox}[Question 006\hfill\emph{solved in Stage 1}]
For a linear Nakayama algebra $A$, denote by $n_A$ the number of indecomposable injective modules over $A$ (up to isomorphism) that have codominant dimension at least 2, and denote by $m_A$ the number of simple modules over $A$ (up to isomorphism) with projective dimension 1 that are not 1-regular. Compute a closed-form formula for $\Sigma_A (n_A - m_A)$, where the sum ranges over all linear Nakayama algebras (up to isomorphism) with $n$ simple modules.
\end{problembox}

\begin{problembox}[Question 007\hfill\emph{solved in Stage 1}]
Let $R$ be the ring of polynomials in $20$ variables and let $I$ be the ideal of $R$ consisting of functions vanishing on the set of points in $\mathbf{C}^{20}$ where some $6$ coordinates are equal to one another. Finally, let $\mathbf{C}$ be the $R$-module which is $\mathbf{C}$ as a vector space, upon which the variables all act by $0$. What is the dimension of $\mathrm{Tor}^2_R(I,\mathbf{C})$?
\end{problembox}

\begin{problembox}[Question 008\hfill\emph{solved in Stage 1}]
For a linear Nakayama algebra $A$ arising from a connected quiver on $n$ vertices, denote by $n_A$ the number of simple (left) modules over $A$ (up to isomorphism) that have projective dimension at least 2. Compute a closed-form formula for $\Sigma_A n_A$, where the sum ranges over all linear Nakayama algebras (up to isomorphism) with $n$ simple modules.
\end{problembox}

\begin{problembox}[Question 009\hfill\emph{solved in Stage 1}]
Let $n\geq 1$ be an integer. For an integer sequence $w=(w_1, \dots , w_{n+1})$ with distinct positive entries and for a permutation $\pi\in S_{n}\,,$ define the integer sequence $a(w,\pi)=(a_1, \dots , a_n)$ as follows. Let 
\[ a_1= \begin{cases}    + \max( w_{\pi(1)}, w_{\pi(1) +1} ) & \text{if } w_{\pi(1)} < w_{\pi(1) +1}\,, \\   - \max( w_{\pi(1)}, w_{\pi(1) +1} ) & \text{if } w_{\pi(1)} > w_{\pi(1) +1}\,,\end{cases}\]
and, if $n\geq 2$, let $(a_2, \dots, a_n) = a(w^-, \pi')$ where $w^-$ is obtained from $w$ by deleting $\max( w_{\pi(1)}, w_{\pi(1) +1} )$, and where $\pi'\in S_{n-1}$ is the permutation defined by
\[ \pi' (i) = \begin{cases}  \pi(i+1) & \text{if } \pi(i+1)  < \pi(1) \,,\\  \pi(i+1)-1 & \text{otherwise} \,.\end{cases}\]
Define the polynomial $F_{w}(y,t)\in \mathbb{Z}[y,t]$ by
\[F_{w}(y,t) = \sum_{\pi \in S_n} y^{\#\{a_i\in a(w,\pi) \mid a_i>0 \} } t^{\#\{i\in\{1,\dots, n-1\} \mid a_i<a_{i+1} \text{ where } a(w,\pi) = (a_1,\dots , a_n) \} } \,.\]
Determine the coefficient of $y^3$ in $F_{w}(y,t)$ for $w=(2, 4, 3, 1, 7, 8, 11, 6, 5, 10, 9)$ as polynomial in $t$.
\end{problembox}

\begin{problembox}[Question 010\hfill\emph{solved in Stage 1}]
Let $c = (1, 2, \dots, n)$ be the standard $n$-cycle in the symmetric group $S_n$. For $p$ satisfying $\gcd(p, n) = 1$ with multiplicative inverse $p'$ so that $pp'=1 \mod n$, let $w \in S_n$ be the permutation with $w(i) \equiv i \cdot p' \pmod n$, so that $w c^p w^{-1} = c.$

Let $B: S_n \to B_n$ be the positive lift of a permutation to the braid group on $n$ strands.  For $n=7$ and $p=5$, express the pure braid $B(w) B(c)^p B(w)^{-1} B(c)^{-1}$ in terms of the Artin generators of the pure braid group $t_{ij} := \sigma_i \sigma_{i+1} \cdots \sigma_{j-2} \sigma_{j-1}^2 \sigma_{j-2}^{-1} \cdots \sigma_{i+1}^{-1} \sigma_i^{-1}$.

For example, if $n=4$ and $p=3$, $w=[3,2,1,4]$ and we should return the list [(2, 3), (1, 3), (1, 2)], representing the pure braid $t_{23} t_{13} t_{12}$.
\end{problembox}

\begin{problembox}[Question 011\hfill\emph{solved in Stage 1}]
Given lattice paths $P$ and $Q$ from $(0, 0)$ to $(m, r)$, with $P$ never rising above $Q$ or incident with $Q$ at any point, we define two matroids in terms of $P$ and $Q$.

For the first matroid, let $P = p_1p_2\cdots p_n$ and $Q = q_1q_2\cdots q_n$. If $p_{u_1}, p_{u_2}, \ldots, p_{u_r}$ are the north steps of $P$ with $u_1 < u_2 < \cdots < u_r$ and $q_{l_1}, q_{l_2}, \ldots, q_{l_r}$ be the north steps of $Q$ with $l_1 < l_2 < \cdots < l_r$.
Let $N_i$ be the interval $[l_i,u_i]$ of integers. Let $M[P,Q]$ be the transversal matroid on the ground set $[n]$ that has the presentation $(N_1,N_2,\ldots,N_r).$ A matroid $M$ for which there exists $P, Q$ such that $M = M[P, Q]$ is called a lattice path matroid.

For the second matroid, label the rows of the diagram enclosed by $P$ and $Q$, from bottom to top, by $1$ through $r$, and the columns, from left to right, by $r+1$ to $r+m$. For each $i \in [r]$, let the set $A_{i}$ consist of $i$ along with the labels of all columns that have a square in row $i$. The rook matroid $R[P, Q]$ is the transversal matroid with the presentation $(A_{1}, \ldots , A_{r})$. A matroid $M$ for which there exists $P, Q$ such that $M = R[P, Q]$ is called a rook matroid.

Let $\lambda/\mu$ be the skew shape enclosed by the region bounded by $P$ and $Q$; then $\lambda/\mu$ determines and is determined by $P, Q$. For $\alpha/\beta = 888888765/76654321$, find the number of lattice path matroids that are isomorphic to the rook matroid defined by $\alpha/\beta$.
\end{problembox}

\begin{problembox}[Question 012\hfill\emph{solved in Stage 1}]
Define a map $s$ from the set of all permutations of length $n$ to that same set as follows.
\begin{enumerate}
\item If $n=1$, then $s(1)=1$.
\item If $p=LnR$, that is $L$ is the substring on the left of the entry $n$, then $s(p)=s(L)s(R)n$.
\end{enumerate}
This operation $s$ is often called the stack sorting.

Let $B_{n,k}$ be the set of permutations $p$ of length $n$ for which all of the following hold.
\begin{enumerate}
\item The permutation $s(p)$ avoids the patterns $2314$ and $3124$.
\item The number of descents of $p$ is $k-1$.
\end{enumerate}

A \emph{block move} is an operation that turns the permutation
\[
p = p_1 p_2 \cdots p_i\ p_{i+1}\cdots p_j\ p_{j+1}\cdots p_k\ p_{k+1}\cdots p_\ell \cdots p_n
\]
into the permutation
\[
p' = p_1 p_2 \cdots p_i\ p_{k+1}\cdots p_\ell\ p_{j+1}\cdots p_k\ p_{i+1}\cdots p_j\ p_{\ell+1}\cdots p_n.
\]
That is, a block move interchanges the string $p_{i+1}\cdots p_j$ and the string $p_{k+1}\cdots p_{\ell}$, for some numbers $i<j<k<\ell$. If $p,p'\in S_n$, and $t$ is the smallest integer so that there exists a sequence of $t$ block moves that turn $p$ into $p'$, then we say that $t=\mathbf{d}(p,p')$ is the \emph{distance} of $p$ and $p'$.

For a given permutation $p$, let $D_{p}(k)$ be the set of all permutations $q$ that are at distance $k$ from $p$.

How many pairs of permutations $(p,\pi)$ are there so that
\begin{enumerate}
\item both $p$ and $\pi$ are of length $30$, and
\item $p\in B_{30,15}$, and
\item $\pi \in D_p(10)$, that is, $\pi$ is at distance $10$ from $p$?
\end{enumerate}
\end{problembox}

\begin{problembox}[Question 013\hfill\emph{solved in Stage 1}]
Fix two positive integers $r$ and $m$. Consider the family of subsets $C_i=\{2i+1,...,2i+r\}$ where $i$ ranges from $0$ to $m-1$. Let $M$ be the matroid of rank $r$ on $n=2(m-1)+r$ elements whose non-bases are exactly the $m$ subsets $C_0,...,C_{m-1}$. Using the notation of the paper \url{https://arxiv.org/pdf/2502.04980} express the unnormalised symmetrized class $[X_{M,sym}]$ in the basis of the symmetric part of the degree $n-r$-part of the Chow ring of the permutohedron given by square free monomials in the differences $H_k:=L_{n-1-k} - L_{n-k}$ (where $k$ runs from $1$ to $n-1$ and we use convention $L_0=0$). As a formula of $r$ and $m$ please give me the coefficient of the monomial $H_{r-1}\cdot \prod_{k=1}^{n-r-1}H_{n-k}$ in this expression.
\end{problembox}

\begin{problembox}[Question 014\hfill\emph{solved in Stage 1}]
Let $\mathcal{A}$ be the hyperplane arrangement that consists of the hyperplanes $\{x_i \pm x_j \mid 1\le i < j \le 11\}\cup \{ x_i \mid 1\le i \le 11 \} $, and let $F_{\mathcal{A}}(y,t)$ be the coarse flag Hilbert-Poincaré series of $\mathcal{A}$.
What is the coefficient of $y^7$ in $F_{\mathcal{A}}(-y,y)$?
\end{problembox}

\begin{problembox}[Question 015\hfill\emph{solved in Stage 1}]
Consider an ordered tree on $n+1$ vertices with root $r$. We impose the breadth-first left-to-right order $<_T$ on the vertices of $T$: vertices are ordered first by distance from the root, and within each level from left to right according to the ordering of children. The covering relation in this order is denoted by $\prec_T$.

Denote by $d(u,v)$ the edge-distance between vertices $u$ and $v$ in $T$. Let $d(T)$ be the maximum of the distances $d(u,v)$ over all pairs $u \prec_T v$ of vertices. Compute a closed-form formula for the number of ordered trees on $n+1$ vertices ($n \geq 10$) such that $d(T) = n - 3$.
\end{problembox}

\begin{problembox}[Question 016\hfill\emph{solved in Stage 1}]
Determine the number of permutations of size $20$ that avoid the patterns $1324$ and $2314$ and have exactly $30$ inversions.
\end{problembox}

\begin{problembox}[Question 017\hfill\emph{solved in Stage 1}]
Fix a multiset
\[
M=\{1^3,2^3,3^4,4,5,6^2,7^2,8^4,9^3,10,11^{13},12^{100},13,14^2,15,16,17,18,19\},
\]
where $i^k$ denotes that $i$ appears with multiplicity $k$.
Let $r=(r_1,\dots,r_{19})$ be a color vector, where each element $i$ is assigned $r_i$ colors.
For each such $r$, consider the colored multiset Eulerian polynomial (i.e., the descent generating polynomial over colored multiset permutations of type $(M,r)$).
For which color vectors $r$ is this polynomial palindromic?
\end{problembox}

\begin{problembox}[Question 018\hfill\emph{solved in Stage 1}]
Let $a(k)$ denote the number of $\{1324, 1342\}$-avoiding permutations of size $k+2$ with exactly $k$ inversions. Determine the generating function $\sum_{k \geq 0} a(k) x^k$.
\end{problembox}

\begin{problembox}[Question 019\hfill\emph{solved in Stage 1}]
Consider a 20 by 20 square grid with lower left corner at (0,0), upper right at (20,20). What is the number of ways to draw 21 monotone grid paths (each having only E or N steps) through the grid which connect (0,0),..,(20,0) with (20,0).....,(20,20) so that no two paths have overlapping steps, endpoints and so that two paths may cross over at a point but not just touch, the path starting at (20,0) ends at (20,12) and the path starting at (10,0) ends at (20,0)?
\end{problembox}

\begin{problembox}[Question 020\hfill\emph{solved in Stage 1}]
For a Dyck path $d$ of semilength $n$, consider its bounce path. The ``reversing bounce path'' of $d$ is the path that arises by reversing $d$, considering its bounce path, and reversing it again. For example, for $d=(()())$, its bounce path is $(())()$ and its reverse bounce path is $()(())$. How many Dyck paths are there for $n=10$, whose bounce path is equal to their reversing bounce path?
\end{problembox}

\begin{problembox}[Question 021\hfill\emph{solved in Stage 1}]
For a permutation $\pi \in S_n$ and $i \in [n]$, denote by $\pi \setminus i$ the permutation that has the same relative order as the sequence obtained by deleting $i$ from $\pi$ written in one-line notation. We call $\pi$ almost decomposable if at least one of the permutations $\pi$, $\pi \setminus 1$, $\pi \setminus \pi_1$, $\pi \setminus n$, or $\pi \setminus \pi_n$ is decomposable.

For all $n \geq 10$, find the smallest value of $k$ such that there exists a $1342$-avoiding permutation of size $n$ with $k$ inversions that is not almost decomposable.
\end{problembox}

\begin{problembox}[Question 022\hfill\emph{solved in Stage 1}]
Let \(Q\) be the poset on $\{1,2,3\}\times\{1,2,3,4\}$ obtained from the product of two chains by deleting exactly the cover relations
\[
(2,1)<(2,2),\qquad (2,3)<(2,4),\qquad (3,2)<(3,3).
\]
For each \(\sigma\in S_4\), label \((p,j)\in Q\) by $\ell_\sigma(p,j)=p+3(\sigma(j)-1)$. If \(\pi=(q_1,\dots,q_{12})\) is a linear extension of \(Q\), define
\[
\operatorname{des}_{\ell_\sigma}(\pi)
=
\#\{i\in\{1,\dots,11\}:\ell_\sigma(q_i)>\ell_\sigma(q_{i+1})\}.
\]
Set
\[
F_Q(x)
=
\sum_{\sigma\in S_4}\sum_{\pi\in L(Q)}
x^{\operatorname{des}_{\ell_\sigma}(\pi)}.
\]

Determine the coefficients \(\gamma_i\) in the expansion
\[
F_Q(x)=\sum_i \gamma_i x^i(1+x)^{\deg F_Q-2i}.
\]
\end{problembox}

\begin{problembox}[Question 023\hfill\emph{solved in Stage 1}]
Let $x$ be a point in the complex Grassmannian $G(k,n)$, and let $M$ be the corresponding matroid. The Chow class of $M$ is defined as the class of the torus orbit closure $\overline{(\mathbb C^*)^n x}$ in the Chow ring of $G(k,n)$.

Expand the Chow class of the lattice path matroid defined by $URUURRUR$ and $RRURRUUU$ as a linear combination of Schubert cycles $s_{\lambda}$.
\end{problembox}

\begin{problembox}[Question 024\hfill\emph{solved in Stage 1}]
Consider the following double cover of graphs: $G_1$ has 5 vertices and edges $\{(1,2),(1,3),(2,3),(2,3),(2,4),(3,5),(4,5),(4,5)\}$, and $G_2$ has $3$ vertices and edges $\{(1,2),(2,2),(2,3),(3,3)\}$. The map on vertices is given by $1 \mapsto 1, 2\mapsto 2, 3 \mapsto 2, 4\mapsto 3, 5\mapsto 3$. Compute the signed graphic matroid of this double cover, and give the set of bases.
\end{problembox}

\begin{problembox}[Question 025\hfill\emph{solved in Stage 1}]
The circle graph with 8 edges has the property that on the Bergman fan of its graphic matroid the coarse fan structure coincides with the minimal nested set structure.

Count up to isomorphism the number of all simple connected graphs $G$ with 8 edges such that the associated graphic matroid $M$ has the property that it is connected and that on the Bergman fan of $M$ the coarse fan structure coincides with the minimal nested set fan structure.
\end{problembox}

\begin{problembox}[Question 026\hfill\emph{solved in Stage 1}]
Let $M$ be the lattice path matroid of the skew shape $(5,5,5,5,5,5,5)/(4,2,2,2,1,1)$.
Compute the Hilbert series of the Chow ring of $M$.
\end{problembox}

\begin{problembox}[Question 027\hfill\emph{solved in Stage 1}]
Determine the normalized volume of the non-nesting rook matroid (see arXiv: 2410.00127) defined by $UUUURURRRRUR$ and $RUURUURRRRUU$.
\end{problembox}

\begin{problembox}[Question 028\hfill\emph{solved in Stage 1}]
Consider the rational normal curve in real projective space $\mathbb P^6$, and consider its image in an affine chart. Place 16 points in convex position on the image of the curve in the affine chart, and take their convex hull. How many distinct combinatorial types of polytopes can you get with this construction?
\end{problembox}

\begin{problembox}[Question 029\hfill\emph{solved in Stage 1}]
Let $\mathcal{A}$ be a hyperplane arrangement of $9$ hyperplanes in a $3$-dimensional complex vector space, with the hyperplanes labeled by $\{0, 1, \dotsc, 8\}$. Suppose that the $1$-dimensional flats of the arrangement are given by the following sets of hyperplanes: $\{0, 1, 2\}, \{0, 3,8\},\{0, 4, 6\},\{0, 5, 7\},\{1, 3, 7\},\{1, 4, 8\},\{1, 5, 6\},\{2, 3, 6\},\{2, 4, 7\},\{2, 5, 8\},\{3, 4, 5\}, \{6, 7, 8\}$. Let $X$ be the visible contours compactification of the projectivization of $\mathcal{A}$, which is naturally embedded in the torus-orbit closure of the point of $Gr(3,9)$ corresponding to $\mathcal{A}$. Let $N$ be the normal bundle of $X$ in this torus-orbit closure. Compute $H^1(X, N)$.
\end{problembox}

\begin{problembox}[Question 030\hfill\emph{solved in Stage 1}]
Consider the log canonical embedding of the moduli space of marked cubic surfaces, that is, marked del Pezzo surfaces of degree 3. What is the dimension of the ambient projective space of this log canonical model? By way of comparison, for the moduli space $M_{0,n}$ that dimension would be $(n-2)!-1$. What is the analogous number for moduli of marked cubic surfaces?
\end{problembox}

\begin{problembox}[Question 031\hfill\emph{solved in Stage 1}]
Let $f: Y \to X$ be a degree 2 cover of smooth, projective, connected curves over the complex numbers. Assume that $L$ is a line bundle on $Y$ such that $f^\ast f_\ast L$ is a stable vector bundle of degree zero on $Y$ and such that $\sigma^\ast L$ is not isomorphic to $L$ for a generator $\sigma$ of the Galois group of the covering $Y/X$.

How many line bundles $M$ on $Y$ have the property that there exists an unramified Galois cover $g: Z \to Y$ of degree three satisfying $g^\ast f^\ast f_\ast L \simeq g^\ast f^\ast f_\ast (L \otimes M)$?
\end{problembox}

\begin{problembox}[Question 032\hfill\emph{solved in Stage 1}]
Consider the tensor train variety parametrized via the map
$$
\Phi:\mathcal P\longrightarrow
V^{=}_{D,d} \subseteq\mathbb{P}\!\left(
\mathbb{C}^{d_1}\otimes\cdots\otimes \mathbb{C}^{d_N}
\right), \quad
(M^1,\dots,M^{N-1},[C])\longmapsto [T]
$$
with
$$
\mathcal P
:=
\prod_{r=1}^{N-1}\mathrm{St}(D_r,n_r)
\times
\mathbb{P}^=\!\left(\mathbb{C}^{D_{N-1}\times d_N}\right)
$$
and $\mathbb{P}^{=}\!\bigl(\mathbb{C}^{D_{N-1}\times d_N}\bigr)$ is the projectivization of the space of full rank $D_{N-1}\times d_N$ complex matrices. Using this parametrization, compute the degree of the tensor train variety with $D = (1, 2, 2, 2, 1)$ and $d = (3, 3, 3, 3)$. Do not use Chern classes.
\end{problembox}

\begin{problembox}[Question 033\hfill\emph{solved in Stage 1}]
Let $\lambda$ be the partition $(2,2,2)$. Let $a=\{t^6,t^5,t^4,t^3,t^2,t\}$ such that, in the notation of this paper \url{https://arxiv.org/abs/0902.1321v3}. Let $T\in \mathrm{SYT}(\lambda/\mu;a)$ be the tableaux with filling $t^6,t^5,t^4,t^3,t^2,t$ read from top to bottom, left to right.

In Theorem 4.2 in this paper \url{https://arxiv.org/abs/0902.1321v3} they say that assuming there exists a filling of the tableaux $(2,2,2)$ with numbers $\omega_1,\ldots,\omega_6$ from left to right and top to bottom such that the determinant of the Jacobian $J$ defined in Theorem 4.2 is non-zero then there exists a linear series in the fiber of $(z+t^6)(z+t^5)(z+t^4)(z+t^3)(z+t^2)(z+t)$ of the Wronski map with certain leading terms of the Pl\"ucker coordinates specified in the theorem if and only if the $\omega_i$ satisfy certain equations.

Can you find numbers $\omega_1,\ldots,\omega_6$ such that the theorem applies that is in the case where there exists a linear series in the fiber of $(z+t^6)(z+t^5)(z+t^4)(z+t^3)(z+t^2)(z+t)$ of the Wronski map. There should be a unique answer in this case. What is the sum $\omega_1+\ldots +\omega_6$?
\end{problembox}

\begin{problembox}[Question 034\hfill\emph{solved in Stage 1}]
Let $E$ be an elliptic curve over the field $\mathbb{C}_p$ without complex multiplication. We assume that $E$ has good supersingular reduction. We define $X$ as the perfectoid cover of $E \times E$ which is the inverse limit along the infinite tower $p: E \times E \rightarrow E \times E$ given by p-multiplication. What is the rank of the N\'eron-Severi group of $X$ over $\mathbb{Z}[1/p]$?
\end{problembox}

\begin{problembox}[Question 035\hfill\emph{solved in Stage 1}]
Let $X$ be the variety of symmetric $6\times 6$ matrices of rank at most $4$, viewed as a subvariety in a projective space $\mathbb{P}^{20}$ of dimension $20$. We are interested in the coisotropic hypersurface $\mathrm{Ch}_3(X)$ which consists of all linear spaces of dimension $5$ in $\mathbb{P}^{20}$ which fail to intersect $X$ transversally. What is the degree of this hypersurface in Pl\"ucker coordinates on the Grassmannian $\mathrm{Gr}(5,\mathbb{P}^{20})$?
\end{problembox}

\begin{problembox}[Question 036\hfill\emph{solved in Stage 1}]
Let $X$ be the Segre embedding of $(\mathbb{P}^1)^{\times 3}$ in $\mathbb{P}^{7}$, and consider its affine cone $CX$ in $\mathbb{C}^8$. Fix coordinates $(x_1,\ldots,x_8)$ for $\mathbb{C}^8$. What is the number of complex critical points of the function $f(x)=2(x_1-u_1)^2+\sum_{i=2}^8(x_i-u_i)^2$ from a generic point $u=(u_1,\ldots,u_8)\in\mathbb{C}^8$ restricted to $CX$?
\end{problembox}

\begin{problembox}[Question 037\hfill\emph{solved in Stage 1}]
Consider the cubic threefold $X$ in $\mathbb{RP}^4$ defined by the following equation:
\[
x_0x_1x_2 + x_1x_2x_3 + x_2x_3x_4 + x_0x_3x_4 + x_0x_1x_4 = 0
\]
Let $Y$ be the union of the two-dimensional linear spaces contained in $X$. How many connected components does $X \setminus Y$ have?
\end{problembox}

\begin{problembox}[Question 038\hfill\emph{solved in Stage 1}]
Let $[L]$ be a general point of the Grassmannian $Gr(6, 17)$, and let $X$ denote the visible contour compactification of the corresponding hyperplane arrangement. Then $X$ is naturally embedded in the torus-orbit closure of $[L]$. Let $E$ denote the conormal bundle of $X$ in this toric variety. Compute the degree of $c_{\dim X}(E)$.
\end{problembox}

\begin{problembox}[Question 039\hfill\emph{solved in Stage 1}]
Let $Z$ be the free graded commutative algebra with one generator in each degree, except in degrees $1$, $4$, and $6$. For a natural number $k$, let $X_k$ be the projective toric variety parametrized by the degree $k$ monomials in $Z$ (inside $\mathbb P(Z_k)$). What is the largest $k$ such that $X_k$ is normal?
\end{problembox}

\begin{problembox}[Question 040\hfill\emph{solved in Stage 1}]
Fix seven general points in the real projective plane. Consider the
21 lines spanned by any two points, the 21 conics spanned by any five
points, and the 7 cubics through all points that pass doubly through
one of the seven points. These $49 = 21+21+7$ curves divide the real
projective plane into regions, called polypols, each of which has a unique canonical differential form. How many of these canonical
forms are linearly independent?
\end{problembox}

\begin{problembox}[Question 041\hfill\emph{solved in Stage 1}]
Let $\Sigma$ be a complete convex fan in $\mathbb{R}^n$ and $X$ its
toric variety with torus $T = (\mathbb{C}^*)^n$ acting on $X$. Fix a point $x^*\in X$ such that the orbit of $x^*$ sits dense in $X$. Let $Z_C \subset
X$ be the torus orbit closure corresponding to the d-dimensional cone $C
\in \Sigma$. We call $X$ equivariantly simple at $Z_C$ iff for every $L
\subset \mathbb{C}^n$ linear space and every sequence $(t_n)_{n \geq 0}
\subset T$ with limit $t=lim_{n\to \infty}(t_n\cdot x^*) \in (Z_C)^{\circ}$ (where $ (Z_C)^{\circ}$ is
the open subset of $(Z_C)$ with no zero coordinates), the limit
$lim_{t\rightarrow \infty} (t_n \cdot L) \in Gr(dim(L),\mathbb{C}^n)$
depends only on $t$, where $t_n \cdot L$ is the inverse coordinate-wise
action of T on L. Let $\Sigma$ be the complete fan normal to the polytope whose vertex coordinates are $(1,2,3),(0,2,3),(1,3,2),(0,3,2),(1,0,2),(1,3,1)$. We realise $X\subseteq \mathbb{P}^7$ as the closure of the image of the map $(\mathbb{C}^*)^3\to \mathbb{P}^7, (t_0,t_1,t_2)\mapsto (t_1^2t_2^3:\:t_1^3t_2^2:\:t_0t_2^2:\:t_0t_1t_2^2:\:t_0t_1^2t_2^2:\:t_0t_1^2t_2^3:\:t_0t_1^3t_2:\:t_0t_1^3t_2^2)$ and let $T=(\mathbb{C}^*)^3$ act on $X$ via the induced action of the multiplication in the domain of this parametrization.
A subdivision of a fan consists of picking a maximal cone of the fan and a hyperplane and replacing the cone by the two cones obtained by intersecting the initial cone with the two half-spaces of the hyperplane.
 What is the minimal number of consecutive subdivisions needed to make $X$ equivariantly simple?
\end{problembox}

\begin{problembox}[Question 042\hfill\emph{solved in Stage 1}]
Let $X\subseteq\mathbb{P}:=\mathbb{P}^3\times\mathbb{P}^3$ be the image of the map $(x,y,z,w) \mapsto \bigl(\,  [xy:z:w:1] \,,\,[zw:x:y:1]\, \bigr)$. In the following, we denote by $\mathbb{G}(k,n)$ the Grassmannian of $k$-dimensional projective subspaces of $\mathbb{P}^n$. Consider the product of Grassmannians $\mathbb{G}:=\mathbb{G}(1,3)\times\mathbb{G}(1,3)$. For a generic $L=L_1\times L_2\in\mathbb{G}$, since $\dim(L)=2=\mathrm{codim}(X)$, we expect that $X\cap L$ is either empty or a $0$-dimensional reduced scheme of $\mathbb{P}$. Let $I_\alpha\subseteq\mathbb{G}$ be the locus of products of linear spaces $L\in\mathbb{G}$ such that the intersection $X\cap L$ is nontransversal (in particular, a $0$-dimensional nonreduced set). The locus $I_\alpha$ is a hypersurface in $\mathbb{G}$, in particular, it has a multidegree $\beta=(\beta_1,\beta_2)$, where $\beta_i$ is the degree in the Pl\"ucker coordinates of the $i$-th Grassmannian $\mathbb{G}(1,3)$. What is the value of $\beta_1+\beta_2$?
\end{problembox}

\begin{problembox}[Question 043\hfill\emph{solved in Stage 1}]
Let $X = \operatorname{Seg}(\mathbb{P}^3 \times \mathbb{P}^3)
\subset \mathbb{P}(\mathbb{C}^4 \otimes \mathbb{C}^4)$ be the Segre variety of rank-1 matrices, $\sigma_2(X)$ its 2nd secant variety (consisting of rank-2 matrices), and
\[
  u = e_1^{\otimes 2} + e_2^{\otimes 2} + \tfrac{1}{100}\,e_3^{\otimes 2}
      + \tfrac{1}{100}\, e_4^{\otimes 2}
  \;\in\; \mathbb{R}^{4 \times 4}.
\]
Consider the polynomial optimization problem
\begin{align}
  \min_{\substack{a^{(k)},\, b^{(k)}\,\in\,\mathbb{R}^4\\k=1,2}}
  &f(a,b)
  \;:=\;
  \sum_{i,j=1}^{4}
  \Bigl(
    u_{ij}
    - \sum_{k=1}^{2} a^{(k)}_i\, b^{(k)}_j
  \Bigr)^{\!2}\\
  \text{subject to } &f(a,b) \le 1/50.
\end{align}
Use the standard Moment-SOS hierarchy for this Polynomial Optimization Problem (POP), where the d-th level optimizes over polynomials of degree at most 2d. If the hierarchy has finite convergence, set s=1. Otherwise, set s=-1. Compute $(s\cdot f_2^*, s\cdot f_3^*, s\cdot f_4^*)$, where $f_{d}^*$ denotes the optimal value of the d-th level of the hierarchy.
\end{problembox}

\begin{problembox}[Question 044\hfill\emph{solved in Stage 1}]
Let $K$ be an algebraically closed and complete extension field of $\mathbb{Q}_p$ and let $X$ be a smooth, proper, connected curve over $K$ of genus $g$. We look at the associated diamond $X^\diamond$ defined by Peter Scholze. Let $\mathcal{O}$ be its structure sheaf. $X^\diamond$ carries an etale topology which we denote by $et$ and a quasi-pro-etale topology which we denote by $qproet$.

Look at the inclusion $Ext^1_{et}(\mathcal{O}, \mathcal{O}) \rightarrow Ext^1_{qproet}(\mathcal{O}, \mathcal{O})$ given by self-extensions of the structure sheaf on $X^\diamond$ with respect to these two topologies. Compute the dimension of the cokernel of this inclusion as a $K$-vector space. If it is not finite dimensional, give back $\infty$.
\end{problembox}

\begin{problembox}[Question 045\hfill\emph{solved in Stage 1}]
Let $N=9$ and $\lambda$ be the partition $(3,2,1)$ and $\mu=(3,1,0)$. Let $a=(t^2,t)$ such that $a^+=a\cup \{0,0,0,0,\infty,\infty,\infty\}$, in the notation of this paper: \href{https://arxiv.org/abs/0902.1321v3}{arXiv:0902.1321v3}. Let $T\in \mathrm{SYT}(\lambda/\mu;a)$ be the tableau with filling of the skew shape with $t$ in the bottom row and $t^2$ in the cell above and to the right of this cell. In the language of the just mentioned paper the tableau is a diagonally increasing skew tableau.

According to Theorem~4.2 in the same paper there exists a linear series $x$, with certain properties. Let $f$ be the map $\mathbb{P}^1\rightarrow \mathbb{P}^2$ corresponding to $x$. Consider $(f,\{t^2,t\})$ as an element in the moduli space of stable maps, with marked points at $t^2$, $t$.

In the limit $t\rightarrow 0$ in the moduli space of stable maps this becomes a degree $4$ plane curve with a degree zero component $X$ attached at $c_\infty$ that has one marked point on it at $c_t$ and one at $c_0$, corresponding to the marked points $t$, $t^2$, respectively. There is a degree $1$ component attached to $X$ at $c_*$.

Assume that $X$ is parametrized such that $c_\infty$ is at $c=\infty$, $c_0$ is at $c=0$, and $c_t$ is at $c=1$. At what value of $c$ is the node $c_*$?
\end{problembox}

\begin{problembox}[Question 046\hfill\emph{solved in Stage 1}]
Imagine you are an observer agent living in a universe with a smooth continuum of visible angular directions. Your visible sky is modeled as the complex projective line $\mathbb P^1(\mathbb C)$, equipped with normalized Fubini-Study probability measure $\sigma$. You as an observer are not told whether this continuum is fundamental, or whether it is the low-energy rendering/simulation of a hidden finite algebraic substrate. If you are a capable observer to figure out whether you live in a simulation or not, you should be able to infer, from a small statistical defect in the visible angular distribution together with a few invariant-theoretic fingerprints, whether the continuum hypothesis is best explained by an underlying quotient substrate. You not given the hidden group, the primitive invariant, or the perturbation strength. These must be reconstructed from the algebraic data and then used in an asymptotic decision problem.

Let $U=\mathbb C^2$, with homogeneous coordinates $(u,v)$, and let $\sigma$ be the normalized Fubini-Study probability measure on \(\mathbb P^1(\mathbb C)\). Under the alternative hypothesis, a hidden computational substrate is assumed to have a non-cyclic finite \(\operatorname{SL}_2(\mathbb C)\)-quotient singularity $X_\Gamma=\operatorname{Spec} \mathbb C[u,v]^\Gamma$, for some unknown finite subgroup $\Gamma\subset \operatorname{SL}_2(\mathbb C).$ The visible angular sky is \(\mathbb P^1(\mathbb C)\), and the local statistical defect is generated by the lowest primitive invariant of the quotient singularity.

The substrate is known only through the following side information.

1. The minimal resolution of $X_\Gamma$ has exceptional intersection graph $E_8$.

2. Equivalently, $\mathbb C[u,v]^\Gamma$ has primitive homogeneous generators of degrees $12,20,30$, satisfying one weighted-homogeneous relation of degree $60$.

3. The degree-$12$ primitive invariant is normalized as follows. It lies in the one-parameter family
\[
F_a(u,v)=u^{11}v+a\,u^6v^6-uv^{11},
\]
and \(a\in \mathbb C\) is chosen so that, if
\[
H_a=\operatorname{Hess}(F_a)
=
\det
\begin{pmatrix}
\dfrac{\partial^2 F_a}{\partial u^2}
&
\dfrac{\partial^2 F_a}{\partial u\partial v}
\\[2mm]
\dfrac{\partial^2 F_a}{\partial v\partial u}
&
\dfrac{\partial^2 F_a}{\partial v^2}
\end{pmatrix},
\]
and
\[
T_a=\operatorname{Jac}(F_a,H_a)
=
\det
\begin{pmatrix}
\dfrac{\partial F_a}{\partial u}
&
\dfrac{\partial F_a}{\partial v}
\\[2mm]
\dfrac{\partial H_a}{\partial u}
&
\dfrac{\partial H_a}{\partial v}
\end{pmatrix},
\]
then the three degree-\(60\) binary forms $F_a^5, H_a^3, T_a^2$, are linearly dependent over $\mathbb C$, with all three coefficients in the linear dependence nonzero. This is the intrinsic syzygy condition seen by the observer.

For such an admissible value of $a$, define the projectively well-defined nonnegative function
\[
Q([u:v])=
\frac{|F_a(u,v)|^2}{\bigl(|u|^2+|v|^2\bigr)^{12}},
\qquad [u:v]\in \mathbb P^1(\mathbb C),
\]
and its mean-zero normalization
\[
\rho([u:v])=
\frac{Q([u:v])}{\displaystyle\int_{\mathbb P^1} Q\,d\sigma}-1.
\]
For sufficiently small \(\varepsilon\), define a probability measure \(P_\varepsilon\) on \(\mathbb P^1(\mathbb C)\) by
\[
\frac{dP_\varepsilon}{d\sigma}=1+\varepsilon\rho.
\]

The null hypothesis is that the visible sky is statistically indistinguishable from the continuum model, $H_0:\varepsilon=0.$ The alternative hypothesis is that the observed sky is the continuum rendering of the hidden quotient substrate, $H_1:\varepsilon=\varepsilon_0>0.$
You as an observer see $n$ independent samples from $P_\varepsilon$, where $\varepsilon\in\{0,\varepsilon_0\}.$ You as an observer scientist use the Neyman/Pearson likelihood-ratio test for
\[
H_0:\varepsilon=0
\qquad\text{against}\qquad
H_1:\varepsilon=\varepsilon_0>0.
\]
In the local regime $\varepsilon_0=\frac{h}{\sqrt n},$ determine the exact constant \(C\) such that, when $h=C,$ the asymptotic Type-I and Type-II errors are both $\alpha=\beta=\Phi(-1).$

Your answer must include the following five items:

1. the hidden group $\Gamma$, up to conjugacy in $\operatorname{SL}_2(\mathbb C)$

2.  the admissible value or values of $a$

3. the exact Fisher information $\mathcal I=\int_{\mathbb P^1}\rho^2\,d\sigma$

4. the exact threshold constant $C=\frac{2}{\sqrt{\mathcal I}}$

5. the limiting log-Bayes factor, in favor of the hidden quotient-substrate sky over the continuum sky, assuming equal prior odds.
\end{problembox}

\begin{problembox}[Question 047\hfill\emph{solved in Stage 1}]
Let $X$ be the variety in $(\mathbb{C}^*)^5$ cut out by $x_3x_4 + x_1 - 1, x_5 + x_2x_3 - 1, x_4x_5 + x_2x_3x_4 + x_1x_2 - 1$. Consider the form
$$ \Omega=
\frac{dx_1 \wedge dx_2}{x_1 x_2}
+
\frac{dx_1 \wedge dx_3}{x_1 x_3}
-
\frac{dx_1 \wedge dx_5}{x_1 x_5}
+
\frac{dx_2 \wedge dx_3}{x_2 x_3}
+
\frac{dx_2 \wedge dx_4}{x_2 x_4}
+
\frac{dx_3 \wedge dx_4}{x_3 x_4}
+
\frac{dx_3 \wedge dx_5}{x_3 x_5}
+
\frac{dx_4 \wedge dx_5}{x_4 x_5}.
$$
Let ${\rm Crit}_X(s)$ be the critical points on $X$ of the multivalued function $\log(x_1^{s_1}\cdots x_5^{s_5})$. Write a closed form in terms of $s$ for
$$ \sum_{x(s) \in {\rm Crit}_X(s)} \Omega(x(s)).$$
\end{problembox}

\begin{problembox}[Question 048\hfill\emph{solved in Stage 1}]
Suppose that I have a graph composed of a chain of 10 cycles of length 20. Each cycle is joined by an edge to the cycle before and after, with the attachment points diametrically opposed to each other. Let $v$ be the vertex on the last cycle diametrically opposed to the attachment point with the previous chain.

Determine the number of degree 4 effective and $v$-reduced divisors $D$ such that $D + 5v$ has Baker--Norine rank at least 3.
\end{problembox}

\begin{problembox}[Question 049\hfill\emph{solved in Stage 1}]
Let $\mathrm{Gr}(2,n)$ be the Grassmannian of lines in $\mathbb{C}^n$ with its Pl\"ucker embedding. Compute the combinatorial rank $\operatorname{cr}(\mathrm{Gr}(2,4),D)$, where $D$ is the union of all Pl\"ucker hyperplane sections.
\end{problembox}

\begin{problembox}[Question 050\hfill\emph{solved in Stage 1}]
Consider the Segre-Veronese embedding of $\mathbb{P}=(\mathbb{P}^1)^{\times 3}$ associated with the line bundle $\mathcal{L}=\mathcal{O}_{\mathbb{P}}(124534,2345729,87685764)$. Let $\mathcal{P}^{100}(\mathcal{L})$ denote the jet bundle of order $100$ of $\mathcal{L}$. For all $i\ge 0$, let $c_i(\mathcal{P}^{100}(\mathcal{L}))$ denote the $i$-th Chern class of $\mathcal{P}^{100}(\mathcal{L})$. Compute the value of
\[
    \sum_{i=0}^3\deg(c_i(\mathcal{P}^{100}(\mathcal{L})))
\]
\end{problembox}

\begin{problembox}[Question 051\hfill\emph{solved in Stage 1}]
Let $f_t:\mathbb{P}^1\rightarrow \mathbb{P}^2$ be a family of maps of degree 3 parametrized by $t$. Consider $(f_t,p_1,p_2,p_3)$ as a family of stable maps with $3$ marked points where the marked points $p_1,p_2,p_3$ are on the flexes of $f_t$. In the case where the $3$ flexes collide as $t\rightarrow 0$ and the domain of the map splits into a degree $2$ component $C$, a degree $1$ component $E_1$, and  a degree $0$ component $E_0$. Assume that $C$ is connected to $E_0$ and $E_1$ is connected to $E_0$. Assume that $E_0$ is parametrized by $c$ such that the node between $C$ and $E_0$ is at $c=\infty$ and there is at least one marked point at $c=0$ and at least one at $c=1$ what is the highest possible value of $c$ of the location of the node between $E_0$ and $E_1$?
\end{problembox}

\begin{problembox}[Question 052\hfill\emph{solved in Stage 1}]
Consider any Dyck path $D$ in the grid $\mathbb{Z}^2$ that stays above the diagonal $(1,0)$ to $(n+1,n)$. Take $P$ as the Minkowski sum of all simplices $\Delta_{ij}={\rm Newt}(1+x_i+x_ix_{i+1}+\ldots+x_i\cdots x_j)$ for all $(i,j)$ on or below $D$. How many facets does $P$ have?
\end{problembox}

\begin{problembox}[Question 053\hfill\emph{solved in Stage 1}]
Suppose I have two functions \( f_1: \mathbb R^d \to \mathbb R^n \), and \( f_2: \mathbb R^n \to \mathbb R \). The functions are of the form \( f_i(x) = \max(A_i x, B_i x) \). Here, \( A_1, B_1 \) are \( n \times d \) matrices with real entries, and \( A_2, B_2 \) are \( 1 \times n \) matrices with non-negative real entries. The composition $f(x) = f_2(f_1(x))$ is convex, continuous, piecewise-linear and positively homogeneous. Its dual polytope is denoted $P_f$. For $d = 2$, $n = 1003$, determine $\max\{ \#\, \text{vertices } P_f : f \text{ is as described above} \}$.
\end{problembox}

\begin{problembox}[Question 054\hfill\emph{solved in Stage 1}]
We consider the space $X$ of all closed $C^2$-submanifolds $S\subset(\mathbb{R}/\mathbb{Z})^2$ of dimension $1$ such that
for every $a\in\mathbb{Z}^2$ with positive entries the map $\pi_a\colon S\to\mathbb{R}/\mathbb{Z},\, x\mapsto\langle x,a\rangle$ is a covering map and the map
        \begin{equation*}
            \mathrm{tr}_a\colon\mathbb{R}/\mathbb{Z}\to(\mathbb{R}/\mathbb{Z})^2,\, y\mapsto\sum_{x\in \pi_a^{-1}(y)}x
        \end{equation*}
is affine linear with slope $(11,5)$. We equip $X$ with the Hausdorff distance. What is the Lebesgue covering dimension of $X$?
\end{problembox}

\begin{problembox}[Question 055\hfill\emph{solved in Stage 1}]
Let $G(a_1,a_2,\ldots,a_n;z)$ be a Goncharov polylogarithm, such that for complex numbers $a_1,a_2,\ldots,a_n$, and $z$, are defined recursively by:
$$
G(;z)=1
$$
$$
G(a_1,a_2,\ldots,a_n;z) = \int_0^z \frac{dt}{t-a_1} G(a_2,\ldots,a_n;t)
$$
with the common shuffle regularization: $G(0;z)=\log(z)$;
Another examples are $G(1;z)=\log(1-z)$, and $G(x;z)=\log(1-z/x)$ where $\log(z)$ is the natural logarithm;

 Now consider the analytic continuation of polylogarithms: $G(x;y) = G(y;x) + i \pi +G(0;y)-G(0;x)$ which is done with a fixed branch choice. Note that this has the same information as the more familiar formula:
$$\log(x-y) = i\pi+\log(y-x)$$

 With that same branch choice, obtain a formula for $G(x,1;y)$ in terms of polylogs where $y$ is only allowed to be on the right of the semicolon when the numbers to the left of the semicolon are 0 or 1, and where $x$ is allowed to be in the right of the semicolon when the numbers to the left of the semicolon are either of $0,1$ or $y$.
\end{problembox}

\begin{problembox}[Question 056\hfill\emph{solved in Stage 1}]
Consider the following matrix:
\[
A =\left(\!\begin{array}{cccccccccccccccc}
     0&1&0&1&0&1&0&1&0&1&0&1&0&1&0&1\\
     0&0&1&1&0&0&1&1&0&0&1&1&0&0&1&1\\
     0&0&0&0&1&1&1&1&0&0&0&0&0&0&0&0\\
     0&0&0&0&0&0&0&0&1&1&1&1&0&0&0&0\\
     0&0&0&0&0&0&0&0&0&0&0&0&1&1&1&1
\end{array}\!\right).
\]
We denote the $i$-th column of $A$ as $\alpha_i$, and we write $A = [\alpha_1,\dots,\alpha_{16}]$. Now, for a fixed choice of $W = (w_1,\dots,w_{16}) \in (\mathbb{C}^*)^{16}$ consider the polynomial map
\[
\Phi_W\colon(\mathbb{C}^*)^{5} \longrightarrow \mathbb{C}^{16}
\]
given by $\Phi_{W}(t) = (w_1 t^{\alpha_1},\dots,w_{16} t^{\alpha_{16}})$. What is the expected maximum likelihood degree of the Zariski closure of the image of this map for the choice $W = (1,1,1,1,1,-1,-1,1,4,-2,-2,1,c_1,c_2,c_3,c_4)$, where $c_1,c_2,c_3,c_4$ are generic complex numbers in the torus?
\end{problembox}

\begin{problembox}[Question 057\hfill\emph{solved in Stage 1}]
What is the maximum likelihood degree of the Grassmannian of lines
 $\mathrm{Gr}(2,n)$ in its Pl\"ucker embedding in $\mathbb P^{\binom{n}{2}-1}$?
\end{problembox}

\begin{problembox}[Question 058\hfill\emph{solved in Stage 1}]
Let \(G_1,\dots,G_N\) be independent standard Gaussian vectors in
\(\mathbb{R}^2\), and set
\[
    S_N:=\frac{1}{N}\sum_{j=1}^N G_jG_j^{\top}.
\]
For \(N\) large enough, define
\[
    Y_j^{(N)}:=S_N^{-1/2}G_j,
    \qquad j=1,\dots,N.
\]
Set
\[
    u_1=(1,0),
    \qquad
    u_2=\left(-\frac{1}{2},\frac{\sqrt3}{2}\right),
    \qquad
    u_3=\left(-\frac{1}{2},-\frac{\sqrt3}{2}\right).
\]
For \(u\in S^1\), define
\[
    Q_N(u)
    :=
    \frac{1}{N}
    \sum_{j=1}^N
    \left|
        \left\langle Y_j^{(N)},u\right\rangle
    \right|^6 .
\]
Finally, let
\[
    B:=
    \begin{pmatrix}
        2 & 1\\
        1 & 3
    \end{pmatrix}.
\]

Determine the exact value of
\[
    \lim_{N\to\infty}
    \frac{1}{N^{1/3}}
    \log
    \mathbb{P}\left(
        Q_N(u_1)\ge 27,\ 
        Q_N(u_2)\ge 27,\ 
        Q_N(u_3)\ge 27,\ 
        N^{1/3}(S_N-I_2)\succeq B
    \right),
\]
where \(\succeq\) denotes the Loewner order on symmetric matrices.
\end{problembox}

\begin{problembox}[Question 059\hfill\emph{solved in Stage 1}]
Consider the definition of Grassmann Distance Degrees of subvarieties of a Grassmannian ${\rm Gr}(2,n)$ as in \href{https://arxiv.org/pdf/2601.22843}{https://arxiv.org/pdf/2601.22843}.

Given a partition $A_1,\dots,A_r$ of $[n]$ with $|A_1| \geq |A_2| \geq \dots \geq |A_r|$, we can define an associated rank-two matroid on $n$ elements by declaring $e\neq f$ to be a non-basis if and only if there exists $j$ such that $e,f \in A_j$. Denote this by $M_{(|A_1|,\dots,|A_r|)}$ and set $m_i = |A_i|$. Such a matroid is always realized by the $2 \times n$ matrix whose columns are partitioned into $r$ blocks of sizes $m_1,\dots,m_r$ (ordered such that the first $m_1$ columns form the first block and so on) and such that two columns are linearly independent if and only if they belong to different blocks. The rowspans of such matrices form a subvariety of the Grassmannian ${\rm Gr}(2,n)$.

What is the Grassmann Distance Degree of the subvariety corresponding to a partition of sizes $(m_1,\dots,m_5) =(2,2,2,2,2)$ in ${\rm Gr}(2,10)$?
\end{problembox}

\begin{problembox}[Question 060\hfill\emph{solved in Stage 2}]
Call a finite dimensional algebra $A$ self-orthogonal if every indecomposable $A$-module $M$ is self-orthogonal meaning $Ext_A^i(M,M)=0$ for all $i>0$. Consider the class of all connected cyclic Nakayama algebras given by quiver and admissible relations over a field $K$ with 8 simple modules that have dominant dimension at least two and are self-orthogonal. What is the sum of the vector space dimensions of those algebras?
\end{problembox}

\begin{problembox}[Question 061\hfill\emph{solved in Stage 2}]
Consider the following symmetric tensor $T\in \mathbb C^{6\times 6\times 6}$:
\begin{align*}
T_{:,:,1} &= \begin{pmatrix} 25 & 20 & 17 & 15 & -16 & -6 \\ 20 & -4 & 2 & -13 & 30 & 27 \\ 17 & 2 & 59 & -28 & 54 & 60 \\ 15 & -13 & -28 & 19 & -2 & 0 \\ -16 & 30 & 54 & -2 & 8 & 15 \\ -6 & 27 & 60 & 0 & 15 & 54 \end{pmatrix}, & T_{:,:,2} &= \begin{pmatrix} 20 & -4 & 2 & -13 & 30 & 27 \\ -4 & 28 & 23 & 35 & -36 & 0 \\ 2 & 23 & 23 & 26 & -18 & -9 \\ -13 & 35 & 26 & 55 & -51 & 3 \\ 30 & -36 & -18 & -51 & 54 & 3 \\ 27 & 0 & -9 & 3 & 3 & -12 \end{pmatrix},\\[1ex]
T_{:,:,3} &= \begin{pmatrix} 17 & 2 & 59 & -28 & 54 & 60 \\ 2 & 23 & 23 & 26 & -18 & -9 \\ 59 & 23 & 27 & 32 & -13 & 18 \\ -28 & 26 & 32 & 62 & -42 & -15 \\ 54 & -18 & -13 & -42 & 49 & 21 \\ 60 & -9 & 18 & -15 & 21 & 12 \end{pmatrix}, & T_{:,:,4} &= \begin{pmatrix} 15 & -13 & -28 & 19 & -2 & 0 \\ -13 & 35 & 26 & 55 & -51 & 3 \\ -28 & 26 & 32 & 62 & -42 & -15 \\ 19 & 55 & 62 & 72 & -43 & 39 \\ -2 & -51 & -42 & -43 & 37 & -27 \\ 0 & 3 & -15 & 39 & -27 & -21 \end{pmatrix},\\[1ex]
T_{:,:,5} &= \begin{pmatrix} -16 & 30 & 54 & -2 & 8 & 15 \\ 30 & -36 & -18 & -51 & 54 & 3 \\ 54 & -18 & -13 & -42 & 49 & 21 \\ -2 & -51 & -42 & -43 & 37 & -27 \\ 8 & 54 & 49 & 37 & -29 & 33 \\ 15 & 3 & 21 & -27 & 33 & 15 \end{pmatrix}, & T_{:,:,6} &= \begin{pmatrix} -6 & 27 & 60 & 0 & 15 & 54 \\ 27 & 0 & -9 & 3 & 3 & -12 \\ 60 & -9 & 18 & -15 & 21 & 12 \\ 0 & 3 & -15 & 39 & -27 & -21 \\ 15 & 3 & 21 & -27 & 33 & 15 \\ 54 & -12 & 12 & -21 & 15 & 9 \end{pmatrix}.
\end{align*}
Determine the rank $r$ and the unique minimum rank decomposition of $T$, so that $T = \sum_{i=1}^{r} a_i^{\otimes 3}$, where $a_1,\ldots,a_r \in \mathbb C^{6}$. As your answer, submit the tuple $(r, \sum_{i=1}^r |a_{i,3}|^2)$, consisting of the rank $r$ and the ``checksum'' of $3$rd coordinates $a_{i,3}$ of the vectors in the minimum rank decomposition.
\end{problembox}

\begin{problembox}[Question 062\hfill\emph{solved in Stage 2}]
Let $p>5$ be a prime and $\mathbb{Z}_p$ the $p$-adic integers. Let $G = {\rm GL}_{12}$ over $\mathbb{Z}_p$ (Dynkin type $A_{11}$). Let $\breve{\mathbb{Z}}_p$ be the integers of the completed unramified closure of $\mathbb{Q}_p$. For $r > 0$ regard $G_r :=G(\breve{\mathbb{Z}}_p/p^r)$ as (the $k:=\overline{\mathbb{F}}_p$ points of) a perfect algebraic group (of perfectly finite type) over the residue field $\mathbb{F}_p$ by using the positive loop functor (also called jet scheme) construction. Let $G_r^1:=ker(G_r \to G_1)$. Denote by $F$ the geometric Frobenius on $G_r$. Let $T$ be an elliptic unramified maximal torus of Coxeter type of $G$ (defined over $\mathbb{Z}_p$) and let $U$ be the unipotent radical of a Borel subgroup containing $T$ (only defined over $\breve{\mathbb{Z}}_p$). We have the corresponding subgroups $T_r^1,U_r^1$ of $G_r^1$, and $T_r^1$ is $F$-stable.

The $k$-variety $X_r = \{g \in G^1_r \colon g^{-1}F(g) \in U^1_r\}$ has an action of $(T^1_r)^F$ by $t : x \mapsto xt$. Let $\phi : (T_r^1)^F \to \Lambda^\times$ be a character, where $\Lambda=\overline{\mathbb{Q}}_\ell$ and $\ell\neq p$ is a prime.

Let $T^1$ be the preimage in $G(\breve{\mathbb{Z}}_p)$ of $T_r^1 \subseteq G_r^1$. Let $\psi$ be any lift to $T(\mathbb{Q}_p)$ of the inflation of $\phi$ to $T^1(\mathbb{Z}_p)$. Assume $\psi$ has a Howe factorization in the sense of Kaletha (see section 3.6 of Kaletha's paper ``Regular supercuspidal representations''), such that $\psi=\psi_{-1}\cdot\psi_0\cdot \psi_1$ with $\psi_i$ of depth $r_i$ with $r_{-1}:=0<r_0<r_1 \leq r$ and the root systems of corresponding twisted Levi sequence is $G^{-1} = T \subsetneq G^0 \subsetneq G^1 = G$, with $G^0$ of (absolute) type $A_5 \times A_5$. (Note that $\psi_{-1}$ will not affect this question, as $\psi|_{T^1(\mathbb{Z}_p)}$ does not depend on the depth zero character $\psi_{-1}$.)

Question: assume $r=20$, $r_0=12$, $r_1=15$. Determine the smallest $i \geq 0$ for which the $\phi$-weight part of the etale cohomology group with compact support $H_c^i(X_r,\Lambda)[\phi]$ is non-zero.
\end{problembox}

\begin{problembox}[Question 063\hfill\emph{solved in Stage 2}]
Let $A$ be the connected linear Nakayama algebra with Kupisch series $[ 5, 5, 4, 4, 3, 3, 3, 2, 1 ]$ and let $T$ denote the trivial extension algebra of $A$ by $D(A)$. What is the sum of all the $\Omega$-periods of the simple $T$-modules?
\end{problembox}

\begin{problembox}[Question 064\hfill\emph{solved in Stage 2}]
For a connected Nakayama algebra $A$ with $n$ simple modules with associated Kupisch series $c(A)=[c_0,\dots,c_{n-1}]$, define $\Sigma A$ as the sum of all entries of the Kupisch series. Let $L_n$ be the set of all connected Nakayama algebras $A$ (linear and cyclic) with $n$ simple modules that are higher Auslander algebras (meaning the global dimension and the dominant dimension of $A$ are finite and coincide) and for which the unique cycle in the resolution quiver of $A$ consists of exactly 1 vertex. Compute $\sum_{A \in L_n}\Sigma A$ for $n=20$.
\end{problembox}

\begin{problembox}[Question 065\hfill\emph{solved in Stage 2}]
Consider the two posets $P$ and $Q$ given by the cover relations
\[
\begin{array}{l@{\qquad}l}
P: & Q:\\[0.3ex]
0\lessdot a_1,a_2,a_3 & 0\lessdot a_1,a_2,a_3\\
a_1\lessdot b_1,b_2 & a_1\lessdot b_1,b_2,b_3\\
a_2\lessdot b_2,b_3 & a_2\lessdot b_2,b_3\\
a_3\lessdot b_3,b_1 & a_3\lessdot b_1\\
b_1\lessdot c_1,c_2 & b_1\lessdot c_1,c_2,c_3\\
b_2\lessdot c_2,c_3 & b_2\lessdot c_2,c_3\\
b_3\lessdot c_3,c_1 & b_3\lessdot c_1\\
c_1\lessdot d_1,d_2 & c_1\lessdot d_1,d_2,d_3\\
c_2\lessdot d_2,d_3 & c_2\lessdot d_2,d_3\\
c_3\lessdot d_3,d_1 & c_3\lessdot d_1,d_3\\
d_1\lessdot e_1,e_2 & d_1\lessdot e_1,e_2,e_3\\
d_2\lessdot e_2,e_3 & d_2\lessdot e_2,e_3\\
d_3\lessdot e_3,e_1 & d_3\lessdot e_1\\
e_1,e_2,e_3 \lessdot 1 & e_1,e_2,e_3 \lessdot 1
\end{array}
\]
Let $\mathit{ex}\Psi_P(y;a,b)$ and $\mathit{ex}\Psi_Q(y;a,b)$ be their respective extended ab-indices, and $\mathit{ex}\Psi_P(1;a,b)$ and $\mathit{ex}\Psi_Q(1;a,b)$ the corresponding pull-back ab-indices.

Determine the smallest coefficient of an ab-monomial in the difference $\mathit{ex}\Psi_P(1;a,b) - \mathit{ex}\Psi_Q(1;a,b)$.
\end{problembox}

\begin{problembox}[Question 066\hfill\emph{solved in Stage 2}]
Consider the cycle graph $C$ on $5$ vertices. Let $L$ be the lattice of flats of the associated matroid and $\mathcal G$ the subset of $L$ given by edge sets of connected induced subgraphs. In the paper \url{https://arxiv.org/abs/2602.21194}, refinements of the doubled nested set fan $\Sigma^{2\mathcal G}$ associated to the pair $(L,\mathcal G)$ are defined. The authors show that the tubing refinement of $\Sigma^{2\mathcal G}$ defines a subdivision of the nested set fan $\Sigma_{\mathcal G}$. Compute the number of maximal cones in this subdivision, using the results from \url{https://arxiv.org/abs/2602.21194}. How many of them are simplicial?
\end{problembox}

\begin{problembox}[Question 067\hfill\emph{solved in Stage 2}]
Consider the definition of the MMC regions and strata in Section 5 of the paper \url{https://arxiv.org/abs/2503.09571}. For $n=15, r=6$ and the sign pattern $\sigma=(+,+,-,+,+,+,+,+,+,+,-,+,+,+,+)$ how many MMC strata $\mathcal{C}_{n,\sigma,r}^0$ of dimension $24$ are there?
\end{problembox}

\begin{problembox}[Question 068\hfill\emph{solved in Stage 2}]
Consider variables $u_{ij}$ labeled by chords in an $n$-gon with vertices $1,\ldots,n$ and the relations $R_{ij}:= u_{ij} + \prod_{kl \nsim ij} u_{kl} -1$ where $kl \nsim ij$ means the chords do not cross. For $n=8$ how many sign patterns in $\{+,-\}^{n(n-3)/2}$ (which assigns a sign to each variable  $u_{ij}$) do not contradict $R_{ij}$?
\end{problembox}

\begin{problembox}[Question 069\hfill\emph{solved in Stage 2}]
Let $ex\Psi_{\Pi_n}(y,a,b)$ be the extended $ab$-index of the set partition lattice $\Pi_n$ of rank $n$.
Compute $ex\Psi_{\Pi_n}(1,1,1)$ for $n=12$.
\end{problembox}

\begin{problembox}[Question 070\hfill\emph{solved in Stage 2}]
Take a square piece of paper of thickness $d$ and side lengths $l=1$ with barycenter at the origin. First, fold the paper randomly and endow the two resulting sides with different topologies. Fold it back to the original structure. Now fold it three times along the longest edge, i.e., half the longest edge. All folding is performed in Euclidean space $\mathbb R^3$.

The structure you have now will have layers on top of each other endowed with different topologies. Now take one of the vertices and fold it to the middle of the longest edge.

Now subdivide the resulting structure into a triangle and a rectangular shape. Puncture the triangle at its barycenter with a pointy steel rod. Unfold the structure again to its original shape with side length 1.

Group the holes according to their distance to the origin. Remember that you have two topologies. This gives the tuple of the cardinality of the distance classes.

Now take the length of this tuple and define it as $c$. Organize all holes on the paper, which has thickness $d$, in a cyclic graph. Let $c$ be the number of possible colors for the edges.
Compute the number of monochromatic perfect matchings of the graph.
\end{problembox}

\begin{problembox}[Question 071\hfill\emph{solved in Stage 2}]
Let $x$ be a point in the complex Grassmannian $G(k,n)$, and let $M$ be the corresponding matroid. The Chow class of $M$ is defined as the class of the torus orbit closure $\overline{(\mathbb C^*)^n x}$ in the Chow ring of $G(k,n)$.

Expand the Chow class of the non-nesting rook matroid (arXiv: 2410.00127) defined by $UUUURURRRRUR$ and $RUURUURRRRUU$ as a linear combination of Schubert cycles $s_{\lambda}$.
\end{problembox}

\begin{problembox}[Question 072\hfill\emph{solved in Stage 2}]
Let $f_1$, ..., $f_{12}$ be 12 linear forms on a 6-dimensional complex vector space V which define a hyperplane arrangement that realizes the matroid FP12. Let W be the subspace of rational functions on PV spanned by elements of the forms $\prod_{i \in B} f_i^{-1}$, where $B$ is a basis of FP12. Choose a generic hyperplane in $V$, and let $X$ be the kernel of the restriction map from $W$ to this generic hyperplane. Compute the dimension of $X$.
\end{problembox}

\begin{problembox}[Question 073\hfill\emph{solved in Stage 2}]
How many combinatorial types of minuscule type $B_3$ strong valuated Coxeter matroids are there? By 'combinatorial type', we mean the number of induced regular subdivisions of the $3$-cube. Give me only the number of such matroids up to symmetry of the cube.
\end{problembox}

\begin{problembox}[Question 074\hfill\emph{solved in Stage 2}]
We use the notation of \url{https://arxiv.org/pdf/2507.02426v2}.

Let $C/K$ be the curve defined by
$$y^2=x(x+1)(x+5)(4x+1)(x+8)(x+8\cdot 2^{1/2})(x+8\cdot2^{3/4})(x+16)(x+32)$$

What are the normalized thicknesses of the double points of the stable marked model of $C$ in the sense of the above article? Return your final answer as a descending list in the form $(a, b, c, \dots)$.
\end{problembox}

\begin{problembox}[Question 075\hfill\emph{solved in Stage 2}]
Let $X$ be the Segre-Veronese embedding of $Y=(\mathbb{P}^1)^6$ with respect to the line bundle $\mathcal{O}_{Y}(1,2,3,4,5,6)$. Compute the Rayleigh-Ritz degree $\mathrm{RRdeg}_{10}(X)$ with respect to a generic quadratic form.
\end{problembox}

\begin{problembox}[Question 076\hfill\emph{solved in Stage 2}]
Let $L$ denote the log-canonical bundle on $\overline{M}_{0,20}$ over a field of characteristic $2$. Compute the dimension of $H^{17}(\overline{M}_{0,20}, L^{-1})$.
\end{problembox}

\begin{problembox}[Question 077\hfill\emph{solved in Stage 2}]
Consider a triangulated surface with $11$ vertices and whose facets are given by the following set of triangles:
\begin{multline*}
\bigl\{(1,5,6),\ (1,2,6),\ (1,5,7),\ (5,10,7),\ (9,10,7),\ (4,7,9),\ (1,7,3),\\
(1,3,4),\ (1,4,8),\ (1,8,9),\ (1,9,10),\ (1,10,2),\ (10,8,2),\ (10,3,8),\\
(10,11,3),\ (2,3,11),\ (2,8,4),\ (4,5,2),\ (6,2,7),\ (2,11,7),\ (8,11,7),\\
(3,7,8),\ (6,7,4),\ (6,4,3),\ (6,3,9),\ (6,9,8),\ (6,8,11),\ (6,11,5),\\
(11,9,5),\ (2,3,5),\ (9,4,11),\ (10,11,4),\ (4,5,10),\ (3,9,5)\bigr\}.
\end{multline*}
To begin, find coordinates for each of the $11$ vertices, such that the topology of the surface is correct (its genus is the right one and facets overlap only on edges or vertices or the empty set; it has to be an embedded simplicial complex) and every triangle is flat (inside a hyperplane and is the convex hull of its vertices). This initial auxiliary answer should be a $3 \times 11$ matrix $A$ of integer numbers. Finally, multiply the $A$ you obtained by $y=[1,1,1]$, $yA=B$ is a $1 \times 11$ list of integer numbers. Call $D$ the number of distinct numbers you get in $B$. Multiply $D$ by the genus of the surface. Print that final number.
\end{problembox}

\begin{problembox}[Question 078\hfill\emph{solved in Stage 2}]
What is the number of equivalence classes of $d$-dimensional virtual polytopes with $d + 1$ vertices up to affine equivalence for $d = 2$?
\end{problembox}

\begin{problembox}[Question 079\hfill\emph{solved in Stage 2}]
Find the number of vertices of the reflexive flow polytope of the complete directed graph with six vertices.
\end{problembox}

\begin{problembox}[Question 080\hfill\emph{solved in Stage 2}]
Let $K$ be the knot $9_{46}$, $n$ an integer and $M_n$ the Dehn filling of $K$ along the slope $1/n$. The two strong inversions of $K$ extend to $M_n$. Let $K^1_n$ and $K^2_n$ be the branch sets in the respective quotients. For $i=1,2$, let $g_i$ be the limit of $Kh(K^i_n)-Kh(K^i_{n-1})$ as $n$ goes to $\infty$, where $Kh(L)$ denotes the (unreduced) Khovanov homology of a link $L$.  What is the dimension of the (hat version of) Heegaard Floer homology of $M_{g_1\cdot g_2}$? All homology theories are taken with coefficients in the field of 2 elements.
\end{problembox}

\begin{problembox}[Question 081\hfill\emph{solved in Stage 2}]
Consider the moduli space of 7 labeled points in linearly general position in the projective plane. This is the open Grassmannian $\mathrm{Gr}(3,7)$ where all 35 Pl\"ucker coordinates are non-zero modulo the action of the 7-dimensional torus. We are interested in the log canonical embedding of this moduli space. What is the dimension of the projective space where this log canonical model lives?
\end{problembox}

\begin{problembox}[Question 082\hfill\emph{solved in Stage 2}]
Let $G$ be an undirected graph with $d$ vertices. Let $\mathcal{I}_G$ be the set of global Markov statements $A\!\perp\!\!\!\!\perp\! B|S$ for $G$ (where ``$\perp\!\!\!\!\perp$'' stands for ``independent'') such that the collection $\{A,B,S\}$ is a partition of $V$.
For any vector $x\in\mathbb{R}^d$ and any subset $S\subseteq[d]$, we denote by $x_S$ the projection of $x$ onto the coordinates with indices in $S$. Fix a subset $\mathcal{X}\subseteq\mathbb{R}^d$ and define $\mathcal{P}^{(0)}_{G}:=\operatorname{conv}(\mathcal{X})\subseteq\mathbb{R}^d$, where ``$\operatorname{conv}$'' stands for convex hull. For $i\ge 1$ define
\[
    \mathcal{P}^{(i)}_{G} := \operatorname{conv}\bigcup_{A\perp\!\!\!\!\perp B|S\in\mathcal{I}_G}\{(x_A,y_B,x_S),(y_A,x_B,x_S)\mid(x_A,x_B,x_S),(y_A,y_B,x_S)\in\mathcal{P}^{(i-1)}_{G}\}\,.
\]
Suppose now that $d=4$ and that $\mathcal{X}$ is the set of $10$ points in $\mathbb{R}^4$ whose coordinates are the columns of the matrix
\[
    \left(\!\begin{array}{cccccccccc}
      \frac{10}{9}&1&1&\frac{4}{5}&\frac{4}{5}&\frac{3}{2}&\frac{5}{9}&\frac{1}{5}&\frac{8}{3}&\frac{2}{5}\\[7pt]
      \frac{3}{5}&\frac{6}{5}&4&\frac{5}{2}&\frac{1}{9}&6&\frac{3}{4}&2&1&\frac{3}{2}\\[7pt]
      \frac{3}{2}&\frac{3}{2}&\frac{3}{5}&3&\frac{4}{3}&\frac{2}{3}&\frac{1}{3}&\frac{4}{5}&2&\frac{1}{2}\\[7pt]
      2&1&\frac{2}{5}&1&\frac{8}{9}&1&1&1&\frac{9}{4}&2
      \end{array}\!\right)\,.
\]
Determine the number of vertices of $\mathcal{P}^{(2)}_{G}$ with respect to the $4$-cycle $G$.
\end{problembox}

\begin{problembox}[Question 083\hfill\emph{solved in Stage 2}]
Consider the $5$-dimensional cube $C_5$ and all the non-simple $4$-dimensional combinatorial types of slices of $C_5$ obtained by intersecting $C_5$ with a hyperplane $H$. How many of these are not combinatorially equivalent to a product of the form $(C_4 \cap H) \times C_1$?
\end{problembox}

\begin{problembox}[Question 084\hfill\emph{solved in Stage 2}]
    Fix \(n\in\mathbb{N}\), \(r\in[n]\), \(q\ge 4\), and
    \(0<t_1<\cdots<t_n\). Set
    \(
    X_n:=2^{\{t_1,\dots,t_n\}}.
    \)
    For \(c=(c_1,\dots,c_p)\in\mathbb{R}^p\), define
    \(\Psi_c:X_n^p\to X_n\) as follows. Given
    \(E=(E_1,\dots,E_p)\), write
    \[
    \bigcup_{j=1}^p E_j=\{s_1<\cdots<s_a\}.
    \]
    Let \(z_0:=0\), and for \(i=1,\dots,a\) put
    \[
    z_i
    :=
    \mathbf{1}_{\{z_{i-1}\le 1\}} \cdot (z_{i-1})_+
    +
    \sum_{j=1}^p c_j\mathbf{1}_{\{s_i\in E_j\}} .
    \]
    Then
    \[
    \Psi_c(E):=\{s_i:1\le i\le a,\ z_i>1\}.
    \]
    
    Let
    \[
    \alpha=(K,n_1,\dots,n_{K-1},B^{(1)},\dots,B^{(K)}),
    \qquad
    n_0=n_K=1,
    \]
    with \(B^{(\ell)}\in\mathbb{R}^{n_\ell\times n_{\ell-1}}\). If
    \(b_i^{(\ell)}\) denotes the \(i\)-th row of \(B^{(\ell)}\), define
    \[
    G_\alpha^{(0)}(E):=E,
    \qquad
    G_\alpha^{(\ell)}(E)
    :=
    \bigl(
    \Psi_{b_1^{(\ell)}}(G_\alpha^{(\ell-1)}(E)),
    \dots,
    \Psi_{b_{n_\ell}^{(\ell)}}(G_\alpha^{(\ell-1)}(E))
    \bigr).
    \]
    Impose
    \(
    \max_{\ell,j}\#\{i:B_{ij}^{(\ell)}\neq 0\}\le q .
    \)
    For \(\rho\in[n]\), let \(\mathcal G_{n,\rho,q}\) be the class of all maps
    \(g:X_n\to X_n\) for which there exists an \(\alpha\) satisfying the above constraint such that $g=G_\alpha^{(K)}$ and $|g(\{t_1,\dots,t_n\})|\le \rho$. Define
    \[
    \mathfrak{N}(n,\rho,q)
    :=
    \sup_{g\in\mathcal G_{n,\rho,q}}
    \min_{\substack{\alpha:\,G_\alpha^{(K)}=g\\
    \max_{\ell,j}\#\{i:B_{ij}^{(\ell)}\neq 0\}\le q}}
    \sum_{\ell=1}^K n_\ell
    \]
    and
    \[
    \mathfrak{S}(n,\rho,q)
    :=
    \sup_{g\in\mathcal G_{n,\rho,q}}
    \min_{\substack{\alpha:\,G_\alpha^{(K)}=g\\
    \max_{\ell,j}\#\{i:B_{ij}^{(\ell)}\neq 0\}\le q}}
    \#\{(\ell,i,j):B_{ij}^{(\ell)}\neq 0\}.
    \]
    Determine the asymptotic order, up to polylogarithmic factors, of
    \[
        \mathfrak{N}(n,n,q)+\mathfrak{S}(n,r,q).
    \]
\end{problembox}

\begin{problembox}[Question 085\hfill\emph{solved in Stage 3}]
Consider the $6$-dimensional linear space of symmetric matrices $L = \operatorname{Span}\{A_1,\dots,A_6\}\subset \mathbb C^{6\times 6}$ generated by
\begin{align*}
A_1 &= \begin{bmatrix}
25&48&-14&18&10&13\\
48&-2&-46&36&-18&-6\\
-14&-46&-2&-6&-5&-16\\
18&36&-6&28&7&12\\
10&-18&-5&7&-22&-1\\
13&-6&-16&12&-1&-11
\end{bmatrix}, & A_2 &= \begin{bmatrix}
48&-2&-46&36&-18&-6\\
-2&89&9&23&47&43\\
-46&9&33&-20&34&12\\
36&23&-20&33&-11&8\\
-18&47&34&-11&19&20\\
-6&43&12&8&20&27
\end{bmatrix},\\[1ex]
A_3 &= \begin{bmatrix}
-14&-46&-2&-6&-5&-16\\
-46&9&33&-20&34&12\\
-2&33&-1&12&25&14\\
-6&-20&12&-26&-21&-6\\
-5&34&25&-21&-1&10\\
-16&12&14&-6&10&18
\end{bmatrix}, & A_4 &= \begin{bmatrix}
18&36&-6&28&7&12\\
36&23&-20&33&-11&8\\
-6&-20&12&-26&-21&-6\\
28&33&-26&27&10&10\\
7&-11&-21&10&6&-2\\
12&8&-6&10&-2&0
\end{bmatrix},\\[1ex]
A_5 &= \begin{bmatrix}
10&-18&-5&7&-22&-1\\
-18&47&34&-11&19&20\\
-5&34&25&-21&-1&10\\
7&-11&-21&10&6&-2\\
-22&19&-1&6&47&11\\
-1&20&10&-2&11&-7
\end{bmatrix}, & A_6 &= \begin{bmatrix}
13&-6&-16&12&-1&-11\\
-6&43&12&8&20&27\\
-16&12&14&-6&10&18\\
12&8&-6&10&-2&0\\
-1&20&10&-2&11&-7\\
-11&27&18&0&-7&62
\end{bmatrix}.
\end{align*}
Let $X \subseteq L$ be the subvariety of matrices of rank $\leq 3$. What is the number of irreducible components of $X$?
\end{problembox}

\begin{problembox}[Question 086\hfill\emph{solved in Stage 3}]
Let $KQ$ be the path algebra where $Q$ is of linearly oriented Dynkin type $A_3$ having 3 vertices. Let $A$ be the total preprojective algebra of $KQ$. How many indecomposable non-zero $A$-modules are there up to isomorphism?
\end{problembox}

\begin{problembox}[Question 087\hfill\emph{solved in Stage 3}]
Let $KQ$ be the path algebra of the Dynkin quiver of type $A_5$ (so $Q$ has 5 vertices) with linear orientation. Let $A$ be the Koszul dual algebra of the Auslander algebra of $KQ$. How many indecomposable spherical $A$-modules are there? Here a module $M$ is called spherical (in the sense of Auslander-Bridger) if $M$ has finite projective dimension equal to $n$ and $Ext_A^i(M,A)=0$ for $i=1,...,n-1$ with the convention that every module of projective dimension less than or equal to 1 is spherical.
\end{problembox}

\begin{problembox}[Question 088\hfill\emph{solved in Stage 3}]
Let $A$ be the connected linear Nakayama algebra with Kupisch series $[2,2,2,2,1]$ given by quiver and relations over a field. How many spherical indecomposable $A$-bimodules are there? Here a module $M$ over an algebra $B$ is called spherical if $M$ has finite projective dimension $n$ and $\mathrm{Ext}_B^i(M,B)=0$ for $i=1,\dots,n-1$. By convention all modules of projective dimension $\leq 1$ is spherical. For this problem, being spherical is meant in the bimodule sense.
\end{problembox}

\begin{problembox}[Question 089\hfill\emph{solved in Stage 3}]
Let $Q$ be the quiver on the vertex set $\{v_1,\dots,v_6\}$ with arrows
\[
\begin{array}{lll}
a_1\colon v_1\to v_2, & a_2\colon v_2\to v_3, & a_3\colon v_2\to v_5,\\
a_4\colon v_3\to v_1, & a_5\colon v_3\to v_4, & a_6\colon v_4\to v_2,\\
a_7\colon v_4\to v_6, & a_8\colon v_5\to v_4, & a_9\colon v_6\to v_3.
\end{array}
\]
Let $A = \mathbb{F}_3 Q / I$, where $\mathbb{F}_3$ is the field with three elements and $I$ is the two-sided ideal of the path algebra $\mathbb{F}_3 Q$ generated by
\[
\begin{array}{l}
a_2 a_4,\quad a_2 a_5 + 2\,a_3 a_8,\quad a_4 a_1 + a_5 a_6,\quad a_6 a_3,\quad a_6 a_2 + 2\,a_7 a_9,\\
a_8 a_6,\quad a_9 a_4,\quad a_9 a_5,\quad a_4 a_1 a_2,\quad a_1 a_2 a_5 a_7.
\end{array}
\]
How many indecomposable $A$-modules are there?
\end{problembox}

\begin{problembox}[Question 090\hfill\emph{solved in Stage 3}]
Let $A$ be a Nakayama algebra with $n$ simples, that is an algebra given by quiver and relations over a field where the quiver is a directed cycle of length $n$ and the relations are monomial relations possibly of length 1, i.e.\ connected cyclic Nakayama algebras and (possibly not connected) linear Nakayama algebras. Determine the number of such Nakayama algebras with 20 simples that have even global dimension and are higher Auslander algebras, i.e.\ where the global dimension and the dominant dimension of the algebra coincide.
\end{problembox}

\begin{problembox}[Question 091\hfill\emph{solved in Stage 3}]
We color each element of the set $\{1,2,\cdots ,1000\}$ red or blue. Then we arrange the red elements in a Standard Young Tableau $S$ and the
blue elements in a Standard Young Tableau $T$ so that the column lengths of both $S$ and $T$ come from the set $\{2,4,6,8,10\}$.
In how many ways is this possible? Provide the first ten digits and the last ten digits of the answer.
\end{problembox}

\begin{problembox}[Question 092\hfill\emph{solved in Stage 3}]
Let $[L]$ be the point of the Grassmannian $Gr(14, 105)$ corresponding to the $14$-dimensional braid arrangement, and let $X$ be the torus-orbit closure of $[L]$. Compute the coefficient of the Schubert variety corresponding to the partition $(91)^{12}, 78, 0$ in the expansion of the homology class of $X$ in the Schubert basis.
\end{problembox}

\begin{problembox}[Question 093\hfill\emph{solved in Stage 3}]
Johnny's a farmer. Johnny has antlers, birds and cats on his farm: At least one cat, more birds than cats and more antlers than birds, and he loves all of them dearly. Johnny likes counting. One day he counts his cats, cubes the number and writes it down in his diary. The next day he counts the birds and the cats. He computes both the sum of cats and birds and the product of cats and birds and multiplies them with each other. The result he writes in his diary.
On the third day, he counts the birds, cubes the number and writes it in his diary.
On the fourth day, he counts the antlers and the cats. He calculates both the sum of cats and antlers and the product of cats and antlers and multiplies them with each other. The result he writes in his diary.
On the fifth day, he counts antlers, birds and cats, multiplies their numbers and writes down the result times 3.
On the sixth day, he counts antlers and birds. He calculates both the sum of birds and antlers and the product of birds and antlers and multiplies them with each other. The result he writes in his diary.
On the seventh day, Johnny only counts the antlers and writes down the cube of their number into his diary. He then sees that it is good, because all seven numbers in his diary add up to 4 times the sum of birds and cats, times the sum of antlers and birds, times the sum of antlers and cats.
How many antlers, birds and cats are on Johnny's farm? Find the solution with the smallest total number of animals. Answer with the 3-tuple (\#antlers, \#birds, \#cats), where \#antlers denotes the number of antlers (likewise for the other animals).
\end{problembox}

\begin{problembox}[Question 094\hfill\emph{solved in Stage 3}]
What is the maximum number of vertices of a hyperplane section of an octahedron? A hyperplane section is the intersection with a hyperplane.
\end{problembox}

\begin{problembox}[Question 095\hfill\emph{solved in Stage 3}]
Let $K_1$ and $K_2$ be any alternating knots with Seifert genus 1 and signature 0. Suppose that $\det(K_1)=2001$ and $\det(K_2)=2025$. Let $M_i$ be the exterior of $K_i$ for $i=1,2$. Let $M$ be the manifold obtained by gluing $M_1$ to $M_2$ along their boundaries, identifying the meridian of one knot with the Seifert longitude of the other knot and vice versa. What is the dimension of the hat version of Heegaard Floer homology of $M$ with coefficients in the field of two elements?
\end{problembox}

\begin{problembox}[Question 096\hfill\emph{solved in Stage 3}]
For $n$ a positive integer define $V(n)$ to be the integer obtained by
using the base 10 digits of $n$ in base 11. I want to evaluate the
series
$\sum_{p \text{ prime}} \frac{1}{V(p)}$ to an accuracy of $10^{-5}$.
\end{problembox}

\begin{problembox}[Question 097\hfill\emph{solved in Stage 3}]
Let \[
I=
\frac{1}{\pi^2}
\int_0^\pi\int_0^\pi
\log \left|1+2e^{ix}+3e^{iy}\right|\,dy\,dx .
\]

Evaluate \(I\) in closed form as a finite expression built from explicit complex numbers, Riemann zeta values, one-variable classical polylogarithm values (\(\operatorname{Li}_s(q)\), where \(s\in\mathbb Z_{>0}\) and \(q\in\mathbb C\)), and the operations \(+,-,\times,\div\).
\end{problembox}

\begin{problembox}[Question 098\hfill\emph{solved in Stage 3}]
In terms of $d$ the dimension, what is the ratio of the surface area of a $d$-dimensional $\ell_p$ ball to its volume?
\end{problembox}

\begin{problembox}[Question 099\hfill\emph{remains unsolved}]
Let $Q$ be the quiver with 4 vertices, with adjacency matrix $\begin{pmatrix} 0 &  2 & 2 & 0 \\ 0 & 0 & 4 & 0\\ 0 & 0 & 0 & 0\\ 2 & 2 & 2 & 0\end{pmatrix}$, and let $d = (1, 3, 4, 1)$.

Count all the unique Luna types for representations of $Q$ of dimension vector $d$, across all stability conditions that admit semistables.
\end{problembox}

\begin{problembox}[Question 100\hfill\emph{remains unsolved}]
Let $P$ be a set of $6$ permutations of size $n$ such that the identity and the reverse permutation are elements of $P$. What is the maximum density of triples shattered by $P$ asymptotically as $n \to \infty$?
\end{problembox}